\documentclass[final]{siamltex}
\usepackage[T1]{fontenc}
\usepackage{ae}
\usepackage{aecompl}
\usepackage[ansinew]{inputenc}

\expandafter\let\csname equation*\endcsname\relax 
\expandafter\let\csname endequation*\endcsname\relax 
\usepackage{amsmath}  
\usepackage{amssymb}
\usepackage{graphicx}
\usepackage{subfig}
\usepackage{mathtools}
\usepackage{epstopdf}
\DeclareGraphicsExtensions{.pdf,.eps,.png,.jpg,.mps}
\usepackage{epsfig}

\newcommand{\rmd}{\mathsf{d}}
\newcommand{\rmi}{\mathsf{i}}

\newcommand{\matr}[1]{\mathsf{#1}}

\newtheorem{defn}{Definition}[section]
\newtheorem{exmp}{\textit{Example}}[section]
\newtheorem{rem}{\textit{Remark}}[section]

\title{Evaluating Matrix Functions by Resummations on Graphs: the Method of Path-Sums}

\author{P-L~Giscard\footnotemark[2] 
\and S~J~Thwaite\footnotemark[2]
\and D~Jaksch\footnotemark[2]\ \footnotemark[3]}

\begin{document}
\maketitle

\renewcommand{\thefootnote}{\fnsymbol{footnote}}

\footnotetext[2]{Clarendon Laboratory, Department of Physics, University of Oxford, Parks
Road, Oxford OX1 3PU, United Kingdom. Email: {\tt p.giscard1@physics.ox.ac.uk}. P-L Giscard is supported by Scatcherd European and EPSRC scholarships. S Thwaite acknowledges support from Balliol College and a Clarendon Scholarship.}
\footnotetext[3]{Centre for Quantum Technologies, National University of Singapore, 3 Science Drive 2, Singapore 117543.}

\renewcommand{\thefootnote}{\arabic{footnote}}

\begin{abstract}
We introduce the method of path-sums which is a tool for analytically evaluating a function of a square discrete matrix, based on the closed-form resummation of infinite families of terms in the corresponding Taylor series. If the matrix is finite, our approach yields the exact result in a finite number of steps. We achieve this by combining a mapping between matrix powers and walks on a weighted directed graph with a universal graph-theoretic result on the structure of such walks. We present path-sum expressions for a matrix raised to a complex power, the matrix exponential, matrix inverse, and matrix logarithm. We present examples of the application of the path-sum method. 
\end{abstract}



\pagestyle{myheadings}
\thispagestyle{plain}
\markboth{P-L~Giscard, S~J~Thwaite and D~Jaksch}{The Method of Path-Sums}

\section{Introduction}
Many problems in applied mathematics, physics, computer science, and engineering are formulated most naturally in terms of matrices, and can be solved  by computing functions of these matrices. Two well-known examples are the use of the matrix inverse in the solution of systems of linear equations, and the application of the matrix exponential to the solution of systems of linear ordinary differential equations with constant coefficients. These applications, among many others, have led to the rise of an active area of research in applied mathematics and numerical analysis focusing on the development of stable and efficient methods for the computation of functions of matrices over $\mathbb{R}$ or $\mathbb{C}$ (see e.g.~\cite{Higham2008}). 

As part of this ongoing effort, we introduce in this article a novel symbolic method for evaluating matrix functions $f$ analytically and in closed form. The method -- which we term the method of path-sums -- is valid for square discrete matrices, and exploits connections between matrix multiplication and graph theory. It is based on the following three central concepts: (i) we describe a method of partitioning a matrix $\matr{M}$ into 
submatrices and associate these with a weighted directed graph $\mathcal{G}$; (ii) we show that the problem of evaluating any submatrix of $f\big(\matr{M}\big)$ is equivalent to summing the weights of all the walks that join a particular pair of vertices in this directed graph; (iii) we use a universal result about the structure of walks on graphs to exactly resum the weights of families of walks to infinite order. This reduces the sum over weighted walks to a sum over weighted paths, a path being forbidden to visit any vertex more than once. For any finite size matrix, the graph $\mathcal{G}$ is finite and so is the number of path.
We apply the method of path-sums to four common matrix functions: a matrix raised to an arbitrary complex power, the matrix inverse, the matrix exponential, and the matrix logarithm. In each case, we obtain an exact closed-form expression that allows the corresponding function $f\big(\matr{M}\big)$ to be analytically evaluated in a finite number of operations, provided $\matr{M}$ has a finite size.

This paper is organized as follows. In \S\ref{sec:RequiredConcepts} we present the foundational material required by the method of path-sums: in \S\ref{sec:I:MatrixPartitions} we describe the partition of a matrix into submatrices 
; in \S\ref{sec:GraphOfMatrixDecomp} we construct the graph corresponding to this partition, and describe the mapping between matrix multiplication and walks on the corresponding graph; in \S\ref{sec:PSresult} we present the closed-form expression for the sum of all the walks on a graph that underpins the method of path-sums. In \S\ref{sec:ComputationOfMatrixFunctions}, we present path-sum expressions for a matrix raised to a complex power, the matrix inverse, the matrix exponential, and the matrix logarithm. These results are proved in Appendix \ref{appendix:MatrixFunctionsProofs}. In \S\ref{sec:examples}, we provide examples of the application of our results. 
\section{Required Concepts}\label{sec:RequiredConcepts}
In this section we present the three main concepts that underpin the method of path-sums. We begin by outlining the partition of an arbitrary matrix $\matr{M}$ into a collection of sub-arrays, and show how this partition leads naturally to the definition of a weighted directed graph $\mathcal{G}$ that encodes the structure of $\matr{M}$. We then show that computing any power of $\matr{M}$ is equivalent to evaluating the weights of a family of walks on $\mathcal{G}$. We conclude by presenting a universal result on the structure of walks on graphs that forms the basis for a closed-form summation of classes of walks on $\mathcal{G}$.

\subsection{\textbf{Matrix partitions}}\label{sec:I:MatrixPartitions}
A partition of a matrix $\matr{M}$ is a regrouping of the elements of $\matr{M}$ into smaller arrays which 
 interpolate between the usual matrix elements of $\matr{M}$ and $\matr{M}$ itself. In this section we
show how these arrays can be used to compute any function that can be expressed as a power series in $\matr{M}$.

\begin{defn}[General matrix partitions]\end{defn}Let $\matr{M}$ be a $D\times D$ matrix over the complex field $\mathbb{C}$. Let $V$ be a $D$-dimensional vector space over $\mathbb{C}$ with orthonormal basis $\{v_i\}$ $(1\le i \le D)$, where we have adopted Dirac notation. Consider an ensemble of vector spaces $V_1,\ldots,V_n$ such that $V_1\oplus V_2\oplus\cdots\oplus V_n = V$. Let $V_j$ have dimension $d_j$ and basis $\{v_{i_{j,k}}\}$, with $1\le k\le d_j$ and $1\le i_{j,k}\le D$, and let $\varepsilon_j$ be the orthogonal projector onto $V_j$, i.e. $\varepsilon_j = \sum_{k=1}^{d_j} v_{i_{j,k}}\,v^{\dagger}_{i_{j,k}}$ where $\dagger$ designates the conjugate transposition.
These projectors satisfy $\varepsilon_i \varepsilon_j = \delta_{i,j}\,\varepsilon_i$, and the closure relation $\sum_{j=1}^n \varepsilon_j = \mathcal{I}$ with $\mathcal{I}$ is the identity operator on $V$. Consider the restriction-operator $\matr{R}_\mu\in\mathbb{C}^{d_\mu\times D}$, such that $\matr{R}_\mu^{\mathrm{T}} \matr{R}_\mu=\varepsilon_\mu$ where $\mathrm{T}$ designates the transposition. 
A general partition of the matrix $\matr{M}$ is then defined to be the ensemble of matrices $\left\{\matr{M}_{\mu\nu}\right\}$ $(1\le (\mu,\nu) \le n)$, where
\begin{align}
  \matr{M}_{\mu\nu} = \matr{R}_\mu \matr{M}\, \matr{R}^{\mathrm{T}}_\nu
\end{align}
is a $d_\mu \times d_\nu$ matrix that defines a linear map $\varphi_{\mu\nu}:V_\nu \rightarrow V_\mu$.
For $\mu\neq \nu$, we call $\matr{M}_{\mu\nu}$ a flip, while for $\mu= \nu$, we call $\matr{M}_{\nu\nu} = \matr{M}_{\nu}$ a static. 
The projectors $\varepsilon_\mu$ are called projector-lattices. In general there is no relationship between $\matr{M}_{\mu\nu}$ and $\matr{M}_{\nu\mu}$. However if $\matr{M}$ is Hermitian then $\matr{M}_{\nu\mu} = \matr{M}_{\mu\nu}^\dagger$ and $\matr{M}_\nu = \matr{M}_\nu^\dagger$. Similar relations can be derived for the case where $\matr{M}$ is symmetric or antisymmetric. 
\begin{rem}[Block matrix representation]\label{rem:block}\end{rem}For any general partition $\{\matr{M}_{\mu\nu}\}$ of $\matr{M}$, there exists a permutation matrix $\matr{P}$ such that all $\matr{M}_{\mu\nu}$ are contiguous blocks in $\matr{P}\matr{M}\matr{P}^{\mathrm{T}}$.

\begin{exmp}[General partition of a matrix]\label{ex:generaldecompo}\end{exmp}To illustrate a general partition, consider the $4 \times 4$ matrix $\matr{M}$ with elements $(\matr{M})_{ij}=m_{ij}$, which can be interpreted as a linear map on the vector space $V =\mathrm{span}\big(v_1,v_2,v_3,v_4\big)$ with $v_1=\left(1,0,0,0\right)^{\mathrm{T}}$,  $v_2=\left(0,1,0,0\right)^{\mathrm{T}}$, etc.
Choosing vector spaces $V_1=\mathrm{span}\big(v_1,v_3,v_4\big)$ and $V_2 = \mathrm{span}\big(v_2\big)$ such that $V_1\oplus V_2 = V$, 
yields the following partition of $\matr{M}$ 
\begin{equation}
\hspace{-.8mm}\matr{M}_{11} \hspace{-.5mm}= \hspace{-1mm}\begin{pmatrix}m_{11}\hspace{-1.2mm}&m_{13}\hspace{-1.2mm}&m_{14}\\m_{31}\hspace{-1.2mm}&m_{33}\hspace{-1.2mm}&m_{34}\\m_{41}\hspace{-1.2mm}&m_{43}\hspace{-1.2mm}&m_{44}\end{pmatrix}\hspace{-.6mm},~\matr{M}_{12} \hspace{-.5mm}=\hspace{-1mm} \begin{pmatrix}m_{12}\\m_{32}\\m_{42}\end{pmatrix}\hspace{-.6mm},~\matr{M}_{21}\hspace{-.5mm}=\hspace{-1mm}\begin{pmatrix}m_{21}&\hspace{-1.2mm}m_{23}&\hspace{-1.2mm}m_{24}\end{pmatrix}\hspace{-.6mm},~\matr{M}_{22}\hspace{-.5mm}=\hspace{-.8mm}\begin{pmatrix}m_{22}\end{pmatrix}\hspace{-.9mm}.\hspace{-4mm}
\end{equation}
\begin{figure}[t!]
\begin{center}
\vspace{-11mm}
\includegraphics[width=.9\textwidth]{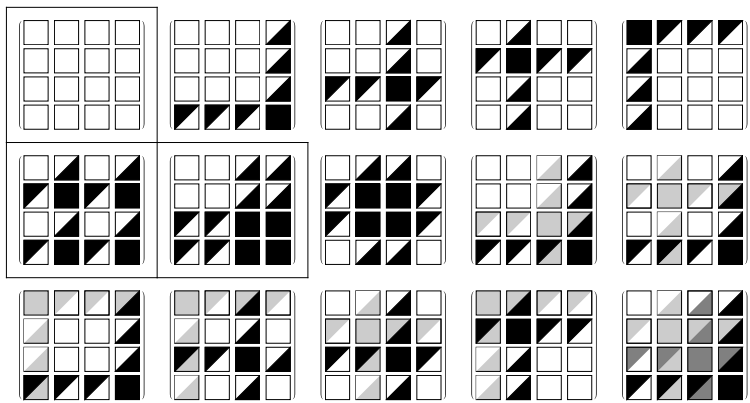}
\vspace{-12mm}
\caption{The $B_4=15$ possible partitions of a $4\times4$ matrix by diagonal projector-lattices. Solid and bicolor squares represent the matrix elements of statics and flips, respectively. The three tensor product partitions are framed. The trivial partition is in the upper-left corner, while the partition into the usual matrix elements is in the lower-right corner.}
\label{fig:MatrixDecomp}
\end{center}
\vspace{-7mm}
\end{figure}
\begin{rem}[Number of general partitions]\end{rem}The partition of a matrix $\matr{M}$ into flips and statics is not unique -- any ensemble of vector spaces such that $\bigoplus_{j=1}^n V_j = V$ produces a valid partition of $\matr{M}$. Consequently a $D\times D$ matrix admits $S(D,n)$ partitions on $n$ such vector spaces, where $S(D,n)$ is the Stirling number of the second kind. It follows that the total number of general partitions of $\matr{M}$ is the $D^{\mathrm{th}}$ Bell number $B_{D}=\sum_{n}S(D,n)$. 
Included in this number are the partition of $\matr{M}$ into the usual matrix elements (i.e.~$\matr{M}_{\mu\nu} = (m_{\mu\nu})$) obtained by choosing $\{V_\mu\}$ to be $D$ subspaces of dimension one, and the partition of $\matr{M}$ into a single static $\matr{M}_{1} = \matr{M}$, obtained by choosing $V_1 = V$. In between these extremes, the subarrays $\matr{M}_{\mu\nu}$ interpolate between the normal matrix elements of $\matr{M}$ and the matrix $\matr{M}$ itself.
Figure \ref{fig:MatrixDecomp} illustrates the whole collection of possible repartitions of a $4\times 4$ matrix into submatrices.

\begin{defn}[Tensor product partitions]\end{defn}
An important subclass of matrix partitions are those that correspond to projector-lattices of tensor product form. This subclass of partitions will be referred to as tensor product partitions, and arises when the vector space $V$ is decomposed as the tensor product (instead of the direct sum) of a collection of subspaces.
Let $V$ be a $D$-dimensional vector space over the complex field $\mathbb{C}$, and consider an ensemble of vector spaces $\mathcal{V}_1,\ldots,\mathcal{V}_N$ such that $\mathcal{V}_1\otimes\cdots\otimes\mathcal{V}_N = V$. Let $\mathcal{V}_i$ have dimension $d_i$ $(2\le d_i \le D)$ and orthonormal basis $\{\mu_{i}\}$ ($1\le \mu \le d_i$), and let $P_{\mu_i}^{(i)} =\mu_{i}\mu^\dagger_{i}$ be the orthogonal projector onto the subspace of $\mathcal{V}_i$ spanned by $\mu_i$. Then a projector-lattice of tensor product form is
\begin{align}
\label{eq:projlatttensorform}
  \varepsilon_\mu^{(S)} = \bigotimes_{i=1}^{S-1} P_{\mu_i}^{(i)}\otimes \mathcal{I}^{(S)}\otimes\bigotimes_{i=S+1}^{N} P_{\mu_i}^{(i)},
\end{align}
where $\mathcal{I}^{(S)}$ is the identity operator on $\mathcal{V}_S$ and $\mu = \left(\mu_1,\ldots,\mu_{S-1},\mu_{S+1},\ldots,\mu_{N}\right)$ is an $(N-1)$-dimensional multi-index denoting which orthogonal projectors are present in $\varepsilon_\mu^{(S)}$. The projector-lattice $\varepsilon_\mu^{(S)}$ acts as a projector on each $\mathcal{V}_i$ $(i\neq S)$ while applying the identity operator to $\mathcal{V}_S$. For fixed $S$ there are $D/d_S$ distinct projector-lattices, corresponding to the different choices of the projector indices $\mu_i$.
For any $D\times D$ matrix $\matr{M}$ and pair of projector-lattices $\varepsilon_\mu^{(S)}$, $\varepsilon_\nu^{(S)}$ there exists a $d_S \times d_S$ matrix $\matr{M}_{\mu\nu}^{(S)}$ such that
\begin{equation}
  \varepsilon_\mu^{(S)}\matr{M}\,\varepsilon_\nu^{(S)} = \bigotimes_{i=1}^{S-1} T^{(i)}_{\mu_i\nu_i}\otimes \matr{M}_{\mu\nu}^{(S)}\otimes\bigotimes_{i=S+1}^N T^{(i)}_{\mu_i\nu_i},
\end{equation}
where $T^{(i)}_{\mu_i\nu_i} = \mu_i\nu^\dagger_i$ is a transition operator from $\nu_i$ to $\mu_i$ in $\mathcal{V}_i$. The matrix $\matr{M}_{\mu\nu}^{(S)}$ defines a linear map on $\mathcal{V}_S$. The ensemble of $(D/d_S)^2$ matrices $\big\{\matr{M}_{\mu\nu}^{(S)}\big\}$ will be referred to as a tensor product partition of $\matr{M}$ on $\mathcal{V}_s$. The three possible tensor product partitions of $\matr{M}$ are illustrated by the framed images in Figure~\ref{fig:MatrixDecomp}.

\subsection{The partition of matrix powers and functions}\label{subsec:PartitionOfMatrixPowers}
Since the matrix elements of $\matr{M}^k$ $(k \in\mathbb{N}^*=\mathbb{N}\backslash \{0\})$ are generated from those of $\matr{M}$ through the rules of matrix multiplication, the partition of a matrix power can be expressed in terms of the partition of the original matrix. Here we present this relationship for the case of a general partition of $\matr{M}$; the case of a tensor product partition is identical. The proof of these results is deferred to Appendix \ref{appendix:MatrixFunctionsProofs}. The partition of $\matr{M}^k$ is given in terms of the partition of $\matr{M}$ by 
$
   \left(\matr{M}^k\right)_{\omega\alpha} = \sum_{\eta_{k},\ldots,\eta_2}\matr{M}_{\omega\eta_{k}} \cdots \matr{M}_{\eta_3\eta_2} \matr{M}_{\eta_2\alpha},
$
where $\alpha \equiv \eta_1$, $\omega \equiv \eta_{k+1}$, and each of the sums runs over the $n$ values that index the vector spaces of the general partition.
It follows that the partition of a matrix function $f(\matr{M})$ with power series expansion $f\big(\matr{M}\big) = \sum_{k=0}^\infty f_k\, \matr{M}^k$ is
\begin{align}
\label{eqn:MatrixFunctionDecomp}
  f\big(\matr{M}\big)_{\omega\alpha} = \sum_{k=0}^\infty f_k\,\sum_{\eta_{k},\ldots,\eta_2}\matr{M}_{\omega\eta_{k}} \cdots \matr{M}_{\eta_3\eta_2} \matr{M}_{\eta_2\alpha}.
\end{align}
This equation provides a method of computing individual submatrices 
of $f\big(\matr{M}\big)$ without evaluating the full result. In the next section, we map the infinite sum of Eq.~\eqref{eqn:MatrixFunctionDecomp} into a sum over the contributions of walks on a weighted graph, thus allowing  exact resummations of families of terms of Eq.~\eqref{eqn:MatrixFunctionDecomp} by applying results from graph theory.

\subsection{\textbf{The graph of a matrix partition}}\label{sec:GraphOfMatrixDecomp}
Given an arbitrary partition of a matrix $\matr{M}$, we construct a weighted directed graph $\mathcal{G}$ that encodes the structure of this partition. Terms that contribute to the matrix power $\matr{M}^k$ are then in one-to-one correspondence with walks of length $k$ on $\mathcal{G}$. The infinite sum over walks on $\mathcal{G}$ involved in the evaluation of $f(\matr{M})$ is then reduced into a sum over paths on $\mathcal{G}$.

\begin{defn}[Graph of a matrix partition]\end{defn}Let $\{\matr{M}_{\mu\nu}\}$ be the partition of $\matr{M}$ formed by a particular set of $n$ projector-lattices $\{\varepsilon_\mu\}$. Then the graph of this matrix partition is defined to be the weighted directed graph $\mathcal{G}=(\mathcal{V},\mathcal{E},\mathrm{w})$, where $\mathcal{V} = \{v_\mu\}$ is a set of $m\leq n$ vertices with the same labels as the projector-lattices, $\mathcal{E} = \{(\nu\mu):\matr{M}_{\mu\nu}\neq \mathbf{0}\}$ is a set of directed edges among these vertices, and $\mathrm{w}$ is an edge-weight function that assigns the submatrix 
$\matr{M}_\mu$ to the loop $(\mu\mu)$ and $\matr{M}_{\mu\nu}$ to the link $(\nu\mu)$. From now on, $\mathcal{G}\backslash\{\alpha,\beta,\ldots\}$ denotes the graph obtained by deleting vertices $\alpha,\beta,\ldots$ from $\mathcal{G}$; and $\mathcal{G}_0$ represents the graph obtained by deleting all self-loops from $\mathcal{G}$.

\begin{rem}[Graph minors]\end{rem}The various graphs that correspond to the different ways of partitioning a matrix $\matr{M}$ into an ensemble of submatrices are minors of the graph obtained by partitioning $\matr{M}$ into its usual matrix elements. Note, this implies that the number of minors obtained by merging vertices\footnote{We allow vertex merging regardless of whether vertices share an edge or not.} on a graph with $D$ vertices is at most $B_D$, with this bound being reached by the complete graph.

\begin{figure}[t!]
\subfloat[]{\includegraphics[scale=0.45]{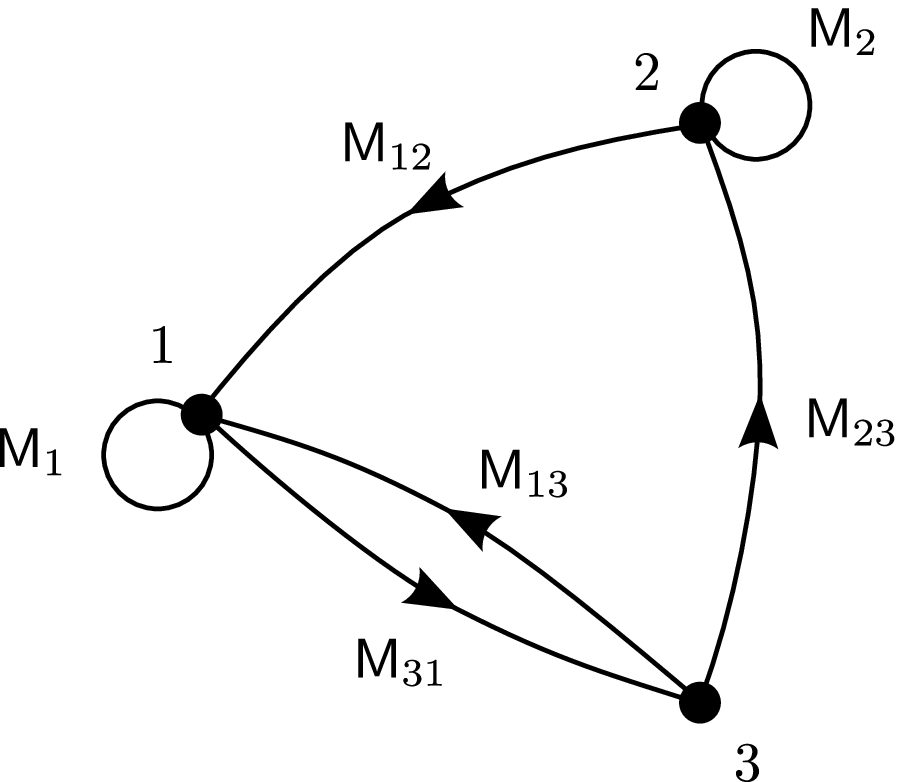}}\hspace{.5cm}
\subfloat[]{\includegraphics[scale=0.45]{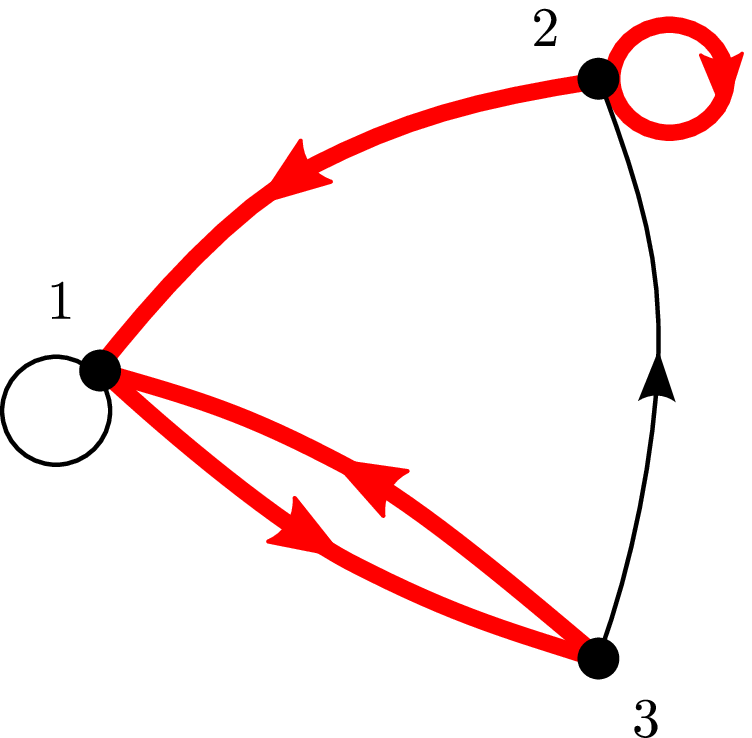}\hspace{.5cm}\includegraphics[scale=0.45]{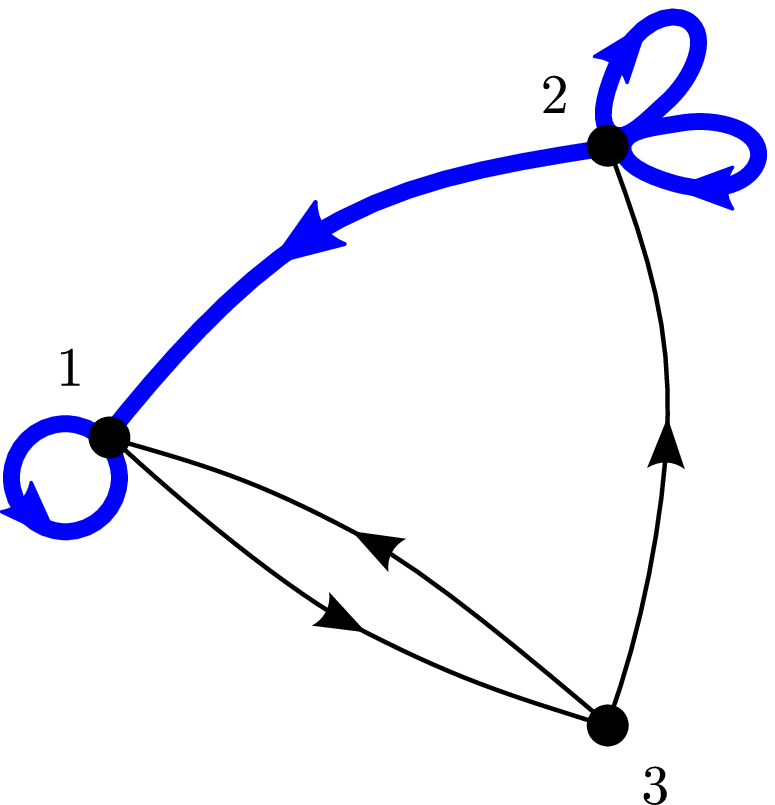}}\hfill
\caption{(a) The graph of the general partition of the $4 \times 4$ matrix in Eq.~\eqref{eqn:MatrixForGraphExample} onto the vector spaces $V_1$, $V_2$, $V_3$ defined in the text. Each edge in the graph is labelled by its weight. (b) Two walks of length 4 from vertex 2 to vertex 1 in thick red (left) and blue (right) lines.
Their contributions $c_{\mathrm{red}}=\matr{M}_{13}\matr{M}_{31}\matr{M}_{12}\matr{M}_{2}$ and $c_{\mathrm{blue}}=\matr{M}_1\matr{M}_{12}\left(\matr{M}_2\right)^2$ form two of the eight terms that sum up to $(\matr{M}^4)_{12}$.}
\label{fig:GraphOfMatrixPartition}
\vspace{-8mm}
\end{figure}
\begin{exmp}[Graph of a matrix partition]\end{exmp}Consider decomposing the $4 \times 4$ matrix
\begin{align}
 \matr{M}=\begin{pmatrix}m_{11}&m_{12}&m_{13}&m_{14}\\m_{21}&m_{22}&m_{23}&m_{24}\\0&0&m_{33}&m_{34}\\m_{41}&m_{42}&0&0\end{pmatrix},
\label{eqn:MatrixForGraphExample}
\end{align}
onto vector spaces $V_1 = \mathrm{span}\big(v_1,v_2)$, $V_2 = \mathrm{span}\big(v_3)$, $V_3 =\mathrm{span}\big(v_4)$ with $v_1=\left(1,0,0,0\right)^{\mathrm{T}}$,  $v_2=\left(0,1,0,0\right)^{\mathrm{T}}$, etc. The corresponding partition of $\matr{M}$ is 
\begin{subequations}
\begin{align}
&\matr{M}_{1}=\begin{pmatrix}m_{11}&m_{12}\\m_{21}&m_{22}\end{pmatrix}, &&\matr{M}_{12}=\begin{pmatrix}m_{13}\\m_{23}\end{pmatrix}, &&\matr{M}_{13}=\begin{pmatrix}m_{14}\\m_{24}\end{pmatrix}\\
&\matr{M}_2 = \begin{pmatrix} m_{33}\end{pmatrix},&&\matr{M}_{23}=\begin{pmatrix}m_{34}\end{pmatrix}, &&\matr{M}_{31}=\begin{pmatrix}m_{41} & m_{42}\end{pmatrix},
\end{align}
\end{subequations}
and $\matr{M}_{21}=\matr{M}_{32}=\matr{M}_{3}=\mathbf{0}$.
Figure \ref{fig:GraphOfMatrixPartition} illustrates $\mathcal{G}$, together with two walks of length 4 from vertex 2 to vertex 1 that contribute to $\big(\matr{M}^4\big)_{12}$.

\begin{defn}[Walks, paths and bare cycles]\end{defn}Consider the graph $\mathcal{G}$ of a matrix partition. Then:\\
\textit{Walk.}~A walk is a vertex-edge sequence, written left-to-right; e.g.~$\left(\alpha\right) \left(\alpha\eta_2\right)\left(\eta_2\right) \cdots\left(\eta_\ell\omega\right)\left(\omega\right)$. A walk that starts and finishes on the same vertex will be termed a closed walk. The set of all walks from $\alpha$ to $\omega$ on $\mathcal{G}$ will be denoted by $W_{\mathcal{G};\alpha\omega}$, and is generally infinite.\\ 
\textit{Path.}~A path is a walk whose vertices are all distinct. The set of all paths from $\alpha$ to $\omega$ on $\mathcal{G}$ will be denoted by $P_{\mathcal{G};\alpha\omega}$. If the graph $\mathcal{G}$ is finite then $P_{\mathcal{G};\alpha\omega}$ is finite.\\
\textit{Bare cycle.}~A bare cycle is a closed walk that does not revisit any internal vertices, i.e. essentially a closed path. The set of all bare cycles off a vertex $\alpha$ on $\mathcal{G}$, denoted by $C_{\mathcal{G};\alpha}$ is finite when $\mathcal{G}$ is finite.

The graph $\mathcal{G}$ provides a useful representation of a matrix partition: each vertex represents a vector space in the partition, while each edge represents a linear mapping between vector spaces. The graph $\mathcal{G}$ is thus a quiver and the ensemble of vector spaces $\{V_\mu\}$ together with the ensemble of linear maps $\{\varphi_{\mu\nu}\}$ is a representation of this quiver. 
Further, the pattern of edges in $\mathcal{G}$ encodes the structure of $\matr{M}$: each loop $(\mu\mu)$ represents a non-zero static, while each link $(\nu\mu)$ represents a non-zero flip. The matrix $\matr{M}$ can therefore be said to be a $\mathcal{G}$-structured matrix. 
Every walk on $\mathcal{G}$, e.g.~$w=(\eta_1)(\eta_1\eta_{2})\cdots (\eta_{k-1}\eta_{k})(\eta_k)$, is now in one-to-one correspondence with a product $c(w)=\matr{M}_{\eta_{k+1}\eta_{k}}\cdots \matr{M}_{\eta_{3}\eta_{2}}\matr{M}_{\eta_{2}\eta_1}$ of submatrices. 
This matrix product is termed the contribution of the walk $w$ and, being a product of matrices, is written right-to-left. This correspondence allows a matrix power to be expressed as a sum over contributions of walks on $\mathcal{G}$. Equation \eqref{eqn:MatrixFunctionDecomp} becomes
\begin{equation}
   f\big(\matr{M}\big)_{\omega\alpha} = \sum_{k=0}^\infty f_k\,\sum_{W_{\mathcal{G};\alpha\omega;k}}c(w),\label{eq:fWalks}
\end{equation}
where $W_{\mathcal{G};\alpha\omega;k}$ is the set of all walks of length $k$ from $\alpha$ to $\omega$ on $\mathcal{G}$.

\subsection{The path-sum}\label{sec:PSresult}
We now state a universal graph theoretic result on the structure of the walks on any graph, which we obtain in \cite{Thwaite2012}. This result reduces a sum of walk contributions, such as the one of Eq.~(\ref{eq:fWalks}), into a sum of weighted paths and bare cycles.
\begin{theorem}[\textbf{Path-sum}]\label{PSresult}
The sum of the contributions of all the walks from $\alpha$ to $\omega$ on a graph $\mathcal{G}$ is given by the path-sum
\begin{subequations}
\begin{align}
 \sum_{w\in W_{\mathcal{G};\alpha\omega}}\hspace{-2mm}c(w)=\sum_{P_{\mathcal{G};\alpha\omega}}c(\omega)'_{\mathcal{G}\backslash\{\alpha,\ldots,\nu_\ell\}}\,c_{\omega\nu_\ell}\cdots c(\nu_2)'_{\mathcal{G}\backslash\{\alpha\}}\,c_{\nu_2\alpha}\, c(\alpha)'_\mathcal{G}\,,\label{eqn:SumOverWalksToPaths2}
\shortintertext{where $c_{\mu\nu}$ is the weight associated to the edge $(\nu\mu)$ and $c(\alpha)'_{\mathcal{G}}$ is given by the sum over the bare cycles}
  c(\alpha)'_\mathcal{G}=\Bigg[\matr{I} -\sum_{C_{\mathcal{G};\alpha}}c_{\alpha\mu_m}\,c(\mu_m)'_{\mathcal{G}\backslash\{\alpha,\ldots,\mu_{m-1}\}}\cdots c_{\mu_3\mu_2} \,
 c(\mu_2)'_{\mathcal{G}\backslash\{\alpha\}} \,c_{\mu_2\alpha}  \Bigg]^{-1},\label{eqn:DressedVertexClosedForm2}
\end{align}
\end{subequations}
with $\matr{I}$ the identity.
\end{theorem}
\noindent The quantity $c(\alpha)'_{\mathcal{G}}$ can be seen as an effective vertex weight resulting from the dressing of vertex $\alpha$ by all the closed walks off $\alpha$ in $\mathcal{G}$.
Remark that if $\mathcal{G}$ is a finite graph, there is only a finite number of paths and bare cycles and Equations (\ref{eqn:SumOverWalksToPaths2}) and (\ref{eqn:DressedVertexClosedForm2}) present finitely many terms. Note that by analytic continuation, Theorem \ref{PSresult} holds regardless of the norm of the walk contributions \cite{Thwaite2012}.

\section{Path-Sum Expressions for Common Matrix Functions}\label{sec:ComputationOfMatrixFunctions}
In \S\ref{sec:RequiredConcepts} we showed that projector-lattices can be used to evaluate the partition of a matrix function $f\big(\matr{M}\big)$, and further, that the resulting expression can be interpreted as a sum over walks on a directed graph $\mathcal{G}$ (Eq.~\eqref{eq:fWalks}). This mapping enables results from graph theory to be applied to the evaluation of matrix power series. In this section we exploit this connection by using Theorem \ref{PSresult} to resum, in closed form, certain families of terms in the power series for the partition of some common matrix functions. For each function we resum all terms in the power series that correspond to closed walks on $\mathcal{G}$, and thereby obtain a closed-form expression for the submatrices $f(\matr{M})_{\omega\alpha}$. Since this expression takes the form of a finite sum over paths on $\mathcal{G}$, we refer to it as a path-sum result. We present path-sum results for a matrix raised to a general complex power, the matrix inverse, matrix exponential, and matrix logarithm. The results are proved in Appendix \ref{appendix:MatrixFunctionsProofs}, and examples illustrating their use are provided in \S\ref{sec:examples}.

\subsection{A matrix raised to a complex power}\label{subsec:MatrixComplexPower}
~\vspace{1mm}
\begin{theorem}[\textbf{Path-sum result for a matrix raised to a complex power}]\label{thm:PathSumMatrixPower}
Let $\matr{M} \in \mathbb{C}^{D\times D}$ be a non-nilpotent matrix, $\{\matr{M}_{\mu\nu}\}$ be an arbitrary partition of $\matr{M}$ and $q \in \mathbb{C}$. Then the partition of $\matr{M}^q$ is given by
\begin{subequations}
\begin{align}
&(\matr{M}^{q})_{\omega\alpha}=-\mathcal{Z}^{-1}\hspace{-1.5mm}\left\{\sum_{P_{\mathcal{G}_0;\alpha\omega}}\matr{F}_{\mathcal{G}\backslash\{\alpha,\ldots,\nu_\ell\}}[\omega]\,{\matr{M}}_{\omega\nu_\ell}\,\cdots \,\matr{F}_{\mathcal{G}\backslash\{\alpha\}}[\nu_2]\,{\matr{M}}_{\nu_2\alpha}\,\matr{F}_\mathcal{G}[\alpha]\right\}\hspace{-1mm}[n]\Bigg|_{n=-q-1}\hspace{-8mm},\label{eqn:MatrixPowerTheorem1} \hspace{-15mm}\\
\shortintertext{where $\mathcal{G}$ is the graph of $\{(\matr{M}-\matr{I})_{\mu\nu}\}$, $\ell$ is the length of the path, and}
&\hspace{-1mm}\matr{F}_{\mathcal{G}}[\alpha] = \left[\matr{I}z^{-1} -\matr{M}_\alpha-\sum_{C_{\mathcal{G}_0;\alpha}}\matr{M}_{\alpha\mu_m} \matr{F}_{\mathcal{G}\backslash{\{\alpha,\ldots,\mu_{m-1}}\}}[\mu_m]\cdots \matr{F}_{\mathcal{G}\backslash\{\alpha\}}[\mu_2]\matr{M}_{\mu_2\alpha}\thinspace\right]^{-1}\hspace{-1mm},\label{eqn:MatrixPowerTheorem2}
\end{align}
with $m$ the length of the bare cycle.
\end{subequations}
\end{theorem}

\noindent Here $\mathcal{Z}^{-1}\left\{g(z)\right\}[n]$ denotes the inverse unilateral Z-transform of $g(z)$. The quantity $\matr{F}_{\mathcal{G}}[\alpha]$ is a matrix continued fraction which terminates at a finite depth if $\mathcal{G}$ is finite.  It is an effective weight associated to vertex $\alpha$ resulting from the dressing of $\alpha$ by all the closed walks off $\alpha$ on $\mathcal{G}$. If $\alpha$ has no neighbours in $\mathcal{G}$ then $\matr{F}_{\mathcal{G}}[\alpha] = [\matr{I}z^{-1}-\matr{M}_{\alpha}]^{-1}$ counts the contributions of all loops on $\alpha$.


\subsection{The matrix inverse}\label{sec:TheMatrixInverse}
~\vspace{1mm}
\begin{theorem}[\textbf{Path-sum result for the matrix inverse}]\label{thm:PathSumMatrixInverse}
Let $\matr{M} \in \mathbb{C}^{D\times D}$ be an invertible matrix, and $\{\matr{M}_{\mu\nu}\}$ be an arbitrary partition of $\matr{M}$. Then as long as all of the required inverses exist, the partition of $\matr{M}^{-1}$ is given by the path-sum
\begin{subequations}
\begin{align}
&\big(\matr{M}^{-1}\big)_{\omega\alpha}=\sum_{P_{\mathcal{G}_0;\alpha\omega}} (-1)^\ell\, \matr{F}_{\mathcal{G}\backslash\{\alpha,\ldots,\nu_\ell\}}[\omega]\matr{M}_{\omega\nu_\ell}\cdots \matr{F}_{\mathcal{G}\backslash\{\alpha\}}[\nu_2]\,\matr{M}_{\nu_2\alpha}\,\matr{F}_{\mathcal{G}}[\alpha],\label{eqn:MatrixInverseTheorem1}\\
\shortintertext{where $\mathcal{G}$ is the graph of $\{(\matr{M}-\matr{I})_{\mu\nu}\}$, $\ell$ is the length of the path, and}
& \matr{F}_{\mathcal{G}}[\alpha] = \left[\matr{M}_{\alpha} - \sum_{C_{\mathcal{G}_0;\alpha}}(-1)^m\,\matr{M}_{\alpha\mu_m} \matr{F}_{\mathcal{G}\backslash{\{\alpha,\ldots,\mu_{m-1}}\}}[\mu_m]\cdots \matr{F}_{\mathcal{G}\backslash\{\alpha\}}[\mu_2]\matr{M}_{\mu_2\alpha}\thinspace\right]^{-1},\label{eqn:MatrixInverseTheorem2}
\end{align}
\end{subequations}
with $m$ the length of the bare cycle.
\end{theorem}

\noindent The quantity $\matr{F}_{\mathcal{G}}[\alpha]$ is a matrix continued fraction which terminates at a finite depth if $\mathcal{G}$ is finite. It is an effective weight associated to vertex $\alpha$ resulting from the dressing of $\alpha$ by all the closed walks off $\alpha$ on $\mathcal{G}$. If $\alpha$ has no neighbours in $\mathcal{G}$ then $\matr{F}_\mathcal{G}[\alpha] = \matr{M}_{\alpha}^{-1}$ counts the contributions of all loops on $\alpha$.

\begin{rem}[Known inversion formulae]\label{rem:KnownFormulae}\end{rem}Two known matrix inversion results can be straightforwardly recovered as special cases of Theorem \ref{thm:PathSumMatrixInverse}. Firstly, by considering the complete directed graph on two vertices, we obtain the well-known block inversion formula
\begin{align}
  \begin{pmatrix}
    \matr{A} & \matr{B} \\ \matr{C} & \matr{D}
  \end{pmatrix}^{-1} = \begin{pmatrix} \left(\matr{A}-\matr{B}\matr{D}^{-1}\matr{C}\right)^{-1} & -\matr{A}^{-1}\matr{B}\left(\matr{D}-\matr{C}\matr{A}^{-1}\matr{B}\right)^{-1} \\ -\matr{D}^{-1}\matr{C}  \left(\matr{A}-\matr{B}\matr{D}^{-1}\matr{C}\right)^{-1} & \left(\matr{D}-\matr{C}\matr{A}^{-1}\matr{B}\right)^{-1}\end{pmatrix}.
\end{align}
Secondly, by applying Theorem \ref{thm:PathSumMatrixInverse} to the linear graph on $N\leq D$ vertices (denoted here $\mathcal{L}_N$), we obtain known continued fraction formulae for the inverse of a $D\times D$ block tridiagonal matrix \cite{Meurant1992, Mallik2001, Kilic2008}. This follows from the observation that a $\mathcal{L}_N$-structured matrix is, up to a permutation of its rows and columns, a block tridiagonal matrix. We provide a general formula for the exponential and inverse of arbitrary $\mathcal{L}_N$-structured matrices in \S\ref{subsec:ExmpLnExp}.

\begin{rem}[Path-sum results via the Cauchy integral formula]\label{rem:CauchyIntegral}\end{rem}A path-sum expression can be derived for any matrix function upon using Theorem \ref{thm:PathSumMatrixInverse} together with the Cauchy integral formula
\begin{align}
   f(\matr{M}) = \frac{1}{2\pi\mathrm{i}}\oint_{\Gamma} f(z)\,(z \matr{I}-\matr{M})^{-1}\,\mathrm{d}z,
\end{align}
where $\rmi^2=-1$, $f$ is an holomorphic function on an open subset $U$ of $\mathbb{C}$ and $\Gamma$ is a closed contour completely contained in $U$ that encloses the eigenvalues of $\matr{M}$. However, for certain matrix functions (including all four we consider in this section) a path-sum expression can be derived independently of the Cauchy integral formula by using Theorem \ref{PSresult} directly on the power series for the function. This method is the one we use to prove the results of this section (see Appendix \ref{appendix:MatrixFunctionsProofs}) and can be extended to matrices over division rings, as will be exposed elsewhere.

\subsection{The matrix exponential}\label{subsec:matrixexp}
~\vspace{1mm}
\begin{theorem}[\textbf{Path-sum result for the matrix exponential}]\label{thm:PathSumMatrixExp}
Let $\matr{M} \in \mathbb{C}^{D\times D}$ and $\{\matr{M}_{\mu\nu}\}$ be an arbitrary partition of $\matr{M}$, with $\mathcal{G}$ the corresponding graph. Then for $\tau \in \mathbb{C}$ the partition of $\exp(\tau\matr{M})$ is given by the path-sum
\begin{subequations}
\begin{align}
 &\hspace{-3mm}\exp(\tau\hspace{0.15mm}\matr{M})_{\omega\alpha} = \mathfrak{L}^{-1}\left\{\sum_{P_{\mathcal{G}_0;\alpha\omega}} \matr{F}_{\mathcal{G}\backslash\{\alpha,\ldots,\nu_\ell\}}[\omega]\matr{M}_{\omega\nu_\ell}\cdots \matr{F}_{\mathcal{G}\backslash\{\alpha\}}[\nu_2]\,\matr{M}_{\nu_2\alpha}\, \matr{F}_\mathcal{G}[\alpha]\right\}(t)\Bigg|_{t\,=\tau}\hspace{-4mm},\hspace{-9mm}\label{eqn:pathsExp}\\
 \shortintertext{where $\ell$ is the length of the path and}
 &\matr{F}_{\mathcal{G}}[\alpha] = \left[s\matr{I}-\matr{M}_{\alpha} -
  \sum_{C_{\mathcal{G}_0;\alpha}}\matr{M}_{\alpha\mu_m} \matr{F}_{\mathcal{G}\backslash{\{\alpha,\ldots,\mu_{m-1}}\}}[\mu_m]\cdots \matr{F}_{\mathcal{G}\backslash\{\alpha\}}[\mu_2]\matr{M}_{\mu_2\alpha}\right]^{-1},\label{eqn:MatrixExponentialTheorem2}
 \end{align}
 with $m$ the length of the bare cycle.
\end{subequations}
\end{theorem}

\noindent Here $s$ is the Laplace variable conjugate to $t$, and $\mathfrak{L}^{-1}\{g(s)\}(t)$ denotes the inverse Laplace transform of $g(s)$. The quantity $\matr{F}_{\mathcal{G}}[\alpha]$ is a matrix continued fraction which terminates at a finite depth if $\mathcal{G}$ is finite. It is an effective weight associated to vertex $\alpha$ resulting from the dressing of $\alpha$ by all the closed walks off $\alpha$ on $\mathcal{G}$. If $\alpha$ has no neighbours in $\mathcal{G}$ then $\matr{F}_{\mathcal{G}}[\alpha] =  [s\matr{I}-\matr{M}_{\alpha}]^{-1}$ counts the contributions of all loops off $\alpha$.
\,\vspace{2mm}
\begin{lemma}[Walk-sum result for the matrix exponential]\label{thm3:ResummedMatrixExponential}
Let $\matr{M} \in \mathbb{C}^{D\times D}$ and $\{\matr{M}_{\mu\nu}\}$ be an arbitrary partition of $\matr{M}$, with $\mathcal{G}$ the corresponding graph. Then for $\tau \in \mathbb{C}$ the partition of $\exp(\tau\matr{M})$ is given by the walk-sum
\begin{multline}\label{eqn:ResummedMatrixExponential}
\exp(\tau\hspace{0.15mm}\matr{M})_{\omega\alpha} =\sum_{W_{\mathcal{G}_0;\alpha\omega}} \int_0^{\tau}\mathrm{d}t_m\cdots\int_0^{t_2}\mathrm{d}t_1\thinspace \exp\big[(t-t_m)\matr{M}_\omega\big]\matr{M}_{\omega\mu_m} \cdots \\\cdots\exp\big[(t_2-t_1)\matr{M}_{\mu_2}\big]\matr{M}_{\mu_2\alpha} \exp\big[ t_1\matr{M}_\alpha\big]\vphantom{\sum_{W_{\mathcal{G}_0;\alpha\omega}} \int_0^{t}}.
\end{multline}
\end{lemma}
\noindent This result corresponds to dressing the vertices only by loops, instead of by all closed walks. An infinite sum over all walks from $\alpha$ to $\omega$ on the loopless graph $\mathcal{G}_0$ therefore remains to be carried out.

\subsection{The matrix logarithm}
~\vspace{1mm}
\begin{theorem}[\textbf{Path-sum result for the principal logarithm}]\label{thm:PathSumMatrixLog}
Let $\matr{M} \in \mathbb{C}^{D\times D}$ be a matrix with no eigenvalues on the negative real axis, and $\{\matr{M}_{\mu\nu}\}$ be a partition of $\matr{M}$. Then as long as all of the required inverses exist, the partition of the principal matrix logarithm of $\matr{M}$ is given by the path-sum
\begin{subequations}
\begin{align}\label{eqn:MatrixLogTheorem1}
  &\big(\log \matr{M}\big)_{\omega\alpha} =\\ &\hspace{2mm}\begin{dcases} \int_0^1\, \mathrm{d}x\,\, x^{-1}\, \left(\matr{I}-\matr{F}_{\mathcal{G}}[\alpha]\right), & \omega = \alpha, \\
  \sum_{P_{\mathcal{G}_0;\alpha\omega}} \int_0^1\, \mathrm{d}x\,\, (-x)^{\ell-1}\, \matr{F}_{\mathcal{G}\backslash\{\alpha,\ldots,\nu_\ell\}}[\omega] \matr{M}_{\omega\nu_\ell} \cdots \matr{F}_{\mathcal{G}\backslash\{\alpha\}}[\nu_2]\, \matr{M}_{\nu_2\alpha}\, \matr{F}_\mathcal{G}[\alpha],& \omega \neq \alpha,
\end{dcases}\nonumber
 \shortintertext{where $\mathcal{G}$ is the graph of $\{(\matr{I}-\matr{M})_{\mu\nu}\}$, $\ell$ the length of the path and}
&\matr{F}_{\mathcal{G}}[\alpha] =\label{eqn:MatrixLogTheorem2}\\
&\hspace{2mm} \left[\matr{I}- x(\matr{I}-\matr{M}_\alpha)-\hspace{-2mm}\sum_{C_{\mathcal{G}_0;\alpha}} (-x)^m\,\matr{M}_{\alpha\mu_m} \matr{F}_{\mathcal{G}\backslash\{\alpha,\ldots,\mu_{m-1}\}}[\mu_m]\,\cdots\, \matr{M}_{\mu_3\mu_2}\matr{F}_{\mathcal{G}\backslash\{\alpha\}} [\mu_2]\matr{M}_{\mu_2\alpha}\right]^{-1}\hspace{-2mm},\nonumber
\end{align}
\end{subequations}
with $m$ the length of the bare cycle.
\end{theorem}

\noindent The quantity $\matr{F}_{\mathcal{G}}[\alpha]$ is a matrix continued fraction which terminates at a finite depth  if $\mathcal{G}$ is finite. It is an effective weight associated to vertex $\alpha$ resulting from the dressing of $\alpha$ by all the closed walks off $\alpha$ on $\mathcal{G}$. If $\alpha$ has no neighbours in $\mathcal{G}$ then $\matr{F}_{\mathcal{G}}[\alpha] = \left[\matr{I}- x(\matr{I}-\matr{M}_\alpha)\right]^{-1}$ counts the contributions of all loops off $\alpha$.

\begin{rem}[Richter relation]\end{rem}The path-sum expression of Theorem \ref{thm:PathSumMatrixLog} is essentially the well-known integral relation for the matrix logarithm \cite{Richter1949, Wouk1965, Higham2008}
 \begin{align}\label{eqn:RichterRelation}
   \log{\matr{M}} = \int_0^1 (\matr{M}-\matr{I}) \big[x(\matr{M}-\matr{I}) +\matr{I}\big]^{-1}\mathrm{d}x,
 \end{align}
with a path-sum expression of the integrand. However, the proof of Theorem \ref{thm:PathSumMatrixLog} that we present in \ref{subsec:MatrixLogProof} does not make explicit use of Eq.~\eqref{eqn:RichterRelation}.\\

\section{Examples}\label{sec:examples}
In this section we present some examples of the application of the path-sum method. In the first section we provide simple numerical examples for a matrix raised to a complex power, the matrix inverse, exponential, and logarithm. In the second part, we provide exact results for the matrix exponential and matrix inverse of block tridiagonal matrices and evaluate the computational cost of path-sum on arbitrary tree-structured matrices.
\subsection{\textbf{Short examples}}

\begin{exmp}[Singular defective matrix raised to an arbitrary complex power]\label{exmp:ExmpMq}\end{exmp}To illustrate the result of Theorem \ref{thm:PathSumMatrixPower}, we consider raising the matrix
\begin{equation}
\matr{M}=\begin{pmatrix} 
-4 & 0 & -1 & 0 & -1 \\
 -2 & -2 & 6 & -2 & 4 \\
 6 & 2 & 1 & -2 & 3 \\
 0 & 0 & -1 & -4 & -1 \\
 -6 & -2 & -5 & 2 & -7 \end{pmatrix},
\end{equation}
to an arbitrary complex power $q$. Note that $\matr{M}$ is both singular and defective; i.e.~non-diagonalizable. We partition $\matr{M}$ onto vector space $V_1 = \mathrm{span}\big(v_1,v_2\big)$ and $V_2 = \mathrm{span}\big(v_3,v_4,v_5\big)$ with $v_1=\left(1,0,0,0\right)^{\mathrm{T}}$ etc., such that $V_1\oplus V_2= V$.
%
The corresponding graph $\mathcal{G}$ is the complete linear graph on two vertices, denoted $\mathcal{K}_2$. Following Theorem \ref{thm:PathSumMatrixPower}, the elements of $\matr{M}^q$ are given by,
\begin{equation}
  (\matr{M}^q)_{ii} = -\mathcal{Z}^{-1}\big\{\matr{F}_{\mathcal{K}_2}[i]\big\}[n]\Bigg|_{n=-q-1}\hspace{-7mm}\textrm{and}~~
  (\matr{M}^q)_{ij} =-\mathcal{Z}^{-1}\big\{ \matr{F}_{\mathcal{K}_2\backslash j}[i]\matr{M}_{12}\matr{F}_{\mathcal{K}_2}[j]  \big\}[n]\Bigg|_{n=-q-1}\hspace{-8mm}, \hspace{0mm}
\end{equation}
where $(i,j)=1,2$, $i\neq j$, $\matr{F}_{\mathcal{K}_2}[i]= \left[\matr{I}z^{-1}-\matr{M}_i-\matr{M}_{ij}\matr{F}_{\mathcal{K}_2\backslash i}[j] \matr{M}_{ji}\right]^{-1}$ and $\matr{F}_{\mathcal{K}_2\backslash i}[j] = \left[\matr{I}z^{-1}-\matr{M}_j\right]^{-1}$.
We thus find in the Z-domain 
\begin{align}
  &\tilde{\matr{M}}(z) =z(4z+1)^{-2}\times\\
   &\begin{pmatrix}
 4 z+1 & 0 & -z & 0 & -z \\
 2 z (4 z-1) & 8 z^2+6 z+1 & \frac{88 z^3+50z^2+6z}{4 z+1} & -2 z (4 z+1) & 
 \frac{56 z^3+34 z^2+4 z}{4 z+1} \\
 2 z (4 z+3) & 2 z (4 z+1) & \frac{88 z^3+54 z^2+13z+1}{4 z+1} & -2 z (4
   z+1) & \frac{56 z^3+22 z^2+3z}{4 z+1} \\
 0 & 0 & -z & 4 z+1 & -z \\
 -2 z (4 z+3) & -2 z (4 z+1) & -\frac{88 z^3+38 z^2+5z}{4 z+1} & 2 z (4
   z+1) & \frac{-56 z^3-6 z^2+5 z+1}{4 z+1}
  \end{pmatrix}.\nonumber
\end{align}
and finally $\matr{M}^q = -\mathcal{Z}^{-1}\{\tilde{\matr{M}}(z)\}[n]|_{n=-q-1}$, which is
\begin{align}\label{eq:MatrixComplexPowerGeneralQ}
\matr{M}^q=\rmi\,(\rmi/2)^{3-2q}\begin{pmatrix}
8 &       0 &      2 q & 0 & 2 q \\
 8 q-4 &    4 &     q(q-2)-11  & 4 & q(q-2)-7  \\
 -8 q-4 &    -4 &   -q(q-2)-3 & 4 & -q(q-2) -7 \\
 0 &         0 &    2 q & 8 & 2 q \\
 8 q+4 & 4 &      q(q-2) +11 & -4 & q(q-2)+15
\end{pmatrix}.
\end{align}
This expression is valid for any $q\in\mathbb{C}$ and fulfills $\matr{M}^{q+q'}=\matr{M}^q\matr{M}^{q'},~\forall(q,q')\in\mathbb{C}^2$. Setting $q=1/2$, we obtain
\begin{equation}
\label{eq:MatrixRoot}
\matr{M}^{1/2}=-\frac{\mathrm{i}}{16}\begin{pmatrix}  32 & 0 & 4 & 0 & 4 \\
 0 & 16 & -47 & 16 & -31 \\
 -32 & -16 & -9 & 16 & -25 \\
 0 & 0 & 4 & 32 & 4 \\
 32 & 16 & 41 & -16 & 57
\end{pmatrix},
\end{equation}
with $\rmi^2=-1$, for which it is easily verified that $(\matr{M}^{1/2})^2=\matr{M}$. Any $p^\mathrm{th}$ root of $\matr{M}$, with $p\in\mathbb{N}^*$, can also be calculated and verified. Further, we note that although $\matr{M}$ is not invertible, setting $q=-1$ in Eq.~\eqref{eq:MatrixComplexPowerGeneralQ} yields the Drazin inverse $\matr{M}^D$ of $\matr{M}$ \cite{Campbell1991}, while setting $q=0$ yields a left and right identity $\matr{M}^\flat$ on $\matr{M}$
\begin{equation}
\matr{M}^{\mathrm{D}}=\frac{1}{16}\begin{pmatrix}
 -4 & 0 & 1 & 0 & 1 \\
 6 & -2 & 4 & -2 & 2 \\
 -2 & 2 & 3 & -2 & 5 \\
 0 & 0 & 1 & -4 & 1 \\
 2 & -2 & -7 & 2 & -9
 \end{pmatrix},\quad
  \matr{M}^\flat = \frac{1}{8}\begin{pmatrix} 8 & 0 & 0 & 0 & 0 \\
 -4 & 4 & -11 & 4 & -7 \\
 -4 & -4 & -3 & 4 & -7 \\
 0 & 0 & 0 & 8 & 0 \\
 4 & 4 & 11 & -4 & 15\end{pmatrix}.
\end{equation}
The above Drazin inverse satisfies indeed $\matr{M}\matr{M}^{\mathrm{D}}\matr{M} = \matr{M}$, $\matr{M}^{\mathrm{D}}\matr{M}\matr{M}^{\mathrm{D}} = \matr{M}^{\mathrm{D}}$ and $\matr{M}^{\mathrm{D}} \matr{M}^q=\matr{M}^q\matr{M}^{\mathrm{D}} =\matr{M}^{q-1}$. We also have $\matr{M}$: $\matr{M}^\flat\matr{M}^q=\matr{M}^q\matr{M}^\flat=\matr{M}^q$ for any $q\in\mathbb{C}$ and finally $\matr{M}^{\mathrm{D}}\matr{M}=\matr{M}\matr{M}^{\mathrm{D}}=\matr{M}^{\flat}$ as expected of the Drazin inverse. These properties imply that $\matr{M}^{\flat}$ is the projector onto $\mathrm{im}(\matr{M})$, the image of $\matr{M}$. 

In addition to the examples with $q = -1$, $0$, and $1/p$ with $p\in\mathbb{N}^*$ presented here, the formula of Eq.~\eqref{eq:MatrixComplexPowerGeneralQ} also holds for any complex value of $q$ and is well behaved. For example, we verify analytically using Theorem \ref{thm:PathSumMatrixPower} to calculate $(\matr{M}^{\pm\rmi})^q$ that $(\matr{M}^{\pm\rmi})^{\mp\rmi}=\matr{M}$.
Finally, it is noteworthy that numerical methods implemented by standard softwares such as MATLAB and \textit{Mathematica} suffer from serious stability problems for the matrix considered here and return incorrect results, as can be seen for the case of $q=1/2$.

\begin{exmp}[Matrix with a non-solvable characteristic polynomial raised to an arbitrary complex power]\label{exmp:GaloisMq}\end{exmp}Consider raising the matrix 
\begin{align}\label{eq:MatrixNonSolvableGalois}
  \matr{M} = \begin{pmatrix}
      -1 & 0 & 0 & 0 & -1/2 \\
 0 & 2 & 5/2 & 2 & -1 \\
 -4 & 0 & 0 & 0 & -1 \\
 0 & -1 & 0 & -2 & 7/4 \\
 0 & 1 & 2 & 0 & 1
  \end{pmatrix},
\end{align}
to an arbitrary complex power.
The characteristic polynomial of $\matr{M}$ is $\chi(x)=x^5-x-1$, whose Galois group is the symmetric group $S_5$ and is thus non-solvable. Proceeding similarly to the example \ref{exmp:ExmpMq} we obtain in the Z-domain
\begin{align}
&\hspace{-3mm}\tilde{\matr{M}}(z)=\frac{z}{z^5+z^4-1}\hspace{-1mm}\left(\begin{matrix}
\hspace{-.7mm} -z (z+1) (z (z+2)-1)-1 &\hspace{-3mm} z^3+\frac{z^2}{2} \\
\hspace{-.7mm} 2 z^2 \left(-4 z^2+z+5\right) &\hspace{-3mm} (z-1) (2 z+1) \left(2 z^2+z+1\right) \\
 \hspace{-.7mm}2 z (z ((z-2) z-2)+2) &\hspace{-3mm} -2 z^4+z^3+z^2 \\
\hspace{-.7mm} \frac{1}{2} z^3 (15 z+8) &\hspace{-3mm} -\frac{1}{4} z (z (z (8 z+3)+7)-4) \\
\hspace{-.7mm} 2 z^2 (z (2 z+5)+4) &\hspace{-3mm} -z (z+1) (2 z+1)
\end{matrix}\right.\hspace{-1mm}\cdots\\
&\hspace{-3mm}\cdots\left.\begin{matrix}
 \frac{1}{4} z^2 (z (2 z+5)+4) &\hspace{-1mm} z^3 &\hspace{-1mm} \frac{1}{2} \left(z-2 z^3\right) \\
 \frac{1}{2} z (z+1)^2 (4 z-5) &\hspace{-1mm} 2 (z-1) z \left(2 z^2+z+1\right) &\hspace{-1mm} -5 z^4+z^3+2 z^2+z \\
 \frac{1}{2} \left(z^2 \left(-z^2+z+4\right)-2\right) &\hspace{-1mm} -2 (z-1) z^3 &\hspace{-1mm} (z-1) z \left(2 z^2-1\right) \\
 -\frac{1}{8} z^2 (z+1) (15 z+8) &\hspace{-1mm} \frac{1}{2} \left(z \left(z \left(-3 z^2+z-4\right)+4\right)-2\right) &\hspace{-1mm} \frac{1}{4} z \left(10 z^3+3 z-7\right) \\
 -\frac{1}{2} z (z+1) (z (2 z+5)+4) &\hspace{-1mm} -2 z^2 (z+1) &\hspace{-1mm} (z+1) \left(2 z^2-1\right)
\end{matrix}\right),\nonumber 
\end{align}
and $\matr{M}^q = -\mathcal{Z}^{-1}\{\tilde{\matr{M}}(z)\}[n]|_{n=-q-1}$. Performing the inverse-Z transform is straightforward and yields $\matr{M}^q=\sum_{k=1}^5 \matr{R}_k(q)$,
with $\matr{R}_k(q)$ a matrix given by
\begin{align}
&\matr{R}_k(q)=\frac{r_k^q}{P(r_k)}\left(\begin{matrix}
 r_k^4-r_k^3+r_k^2+3 r_k+1 & -\frac{1}{2} r_k (r_k+2) \\
 -2 (r_k+1) (5 r_k-4) & (r_k-1) (r_k+2) \left(r_k^2+r_k+2\right) \\
 -2 \left(2 r_k^3-2 r_k^2-2 r_k+1\right) & (1-r_k) (r_k+2) \\
 \frac{1}{2} (-8 r_k-15) & \frac{1}{4} \left(-4 r_k^3+7 r_k^2+3 r_k+8\right) \\
 -2 \left(4 r_k^2+5 r_k+2\right) & r_k (r_k+1) (r_k+2)
\end{matrix}\right.\cdots\\
&\hspace{-2mm}\cdots\left.
\begin{matrix} 
\frac{1}{4} \left(-4 r_k^2-5 r_k-2\right) &\hspace{-1mm} -r_k &\hspace{-1mm} -\frac{1}{2} r_k \left(r_k^2-2\right) \\
 \frac{1}{2} (r_k+1)^2 (5 r_k-4) & -2 r_k \left(r_k^2+1\right) \left(2 r_k^2-3\right) &\hspace{-1mm} -r_k^3-2 r_k^2-r_k+5 \\
 \frac{1}{2} (r_k+1) \left(2 r_k^3-2 r_k^2-2 r_k+1\right) &\hspace{-1mm} -2 (r_k-1) &\hspace{-1mm} (1-r_k) \left(r_k^2-2\right) \\
 \frac{1}{8} (r_k+1) (8 r_k+15) &\hspace{-1mm} \frac{1}{2} \left(2 r_k^4-4 r_k^3+4 r_k^2-r_k+3\right) &\hspace{-1mm} \frac{1}{4} \left(7 r_k^3-3
   r_k^2-10\right) \\
 \frac{1}{2} (r_k+1) \left(4 r_k^2+5 r_k+2\right) &\hspace{-1mm} 2 r_k (r_k+1) &\hspace{-1mm} r_k (r_k+1) \left(r_k^2-2\right)
 \end{matrix}\right)\nonumber,
\end{align}
where $r_k$ is the $k$th root of $\chi(x)$ and $P(r_k)=\prod_{i=1, i\neq k}^5(r_k-r_i)$. Using the analytical properties of the roots of $\chi(x)$, we verify that the above formula yields the correct integer-powers, inverse and $p$th-roots ($p\in\mathbb{N}^*$) of $\matr{M}$.

\begin{exmp}[Matrix inverse]\end{exmp}To illustrate the application of Theorem \ref{thm:PathSumMatrixInverse} we compute the inverse of the matrix of Eq.\eqref{eq:MatrixNonSolvableGalois}.
We partition $\matr{M}$ onto vector spaces $V_1=\mathrm{span}\big(v_1\big)$, $V_2 = \mathrm{span}\big(v_3,v_5\big)$ and $V_3 = \mathrm{span}\big(v_2,v_4\big)$, giving
\begin{align}
&\matr{M}_{1}=\begin{pmatrix}-1\end{pmatrix},\hspace{4.4mm}\matr{M}_{2}=\begin{pmatrix}0 & -1 \\
 2 & 1\end{pmatrix}, \hspace{6mm}\matr{M}_{3}=\begin{pmatrix} 2 & 2 \\
 -1 & -2\end{pmatrix},\\
&\matr{M}_{12}=\begin{pmatrix}0 & -1/2\end{pmatrix},~\matr{M}_{21}=\begin{pmatrix}-4 \\
 0\end{pmatrix},~\matr{M}_{23}=\begin{pmatrix} 0 & 0 \\
 1 & 0\end{pmatrix},~\matr{M}_{32}=\begin{pmatrix}5/2 & -1 \\
 0 & 7/4\end{pmatrix},
\end{align}
and $\matr{M}_{13}=\matr{M}_{31}=0$.
The corresponding graph is thus the linear graph on three vertices, denoted by $\mathcal{L}_3$. By Theorem \ref{thm:PathSumMatrixInverse}, the diagonal elements of $\matr{M}^{-1}$ are given by
\begin{subequations}
\begin{align}
  \left(\matr{M}^{-1}\right)_{11} = \matr{F}_{\mathcal{L}_3}[1], \qquad \left(\matr{M}^{-1}\right)_{22} = \matr{F}_{\mathcal{L}_3}[2], \qquad \left(\matr{M}^{-1}\right)_{33} = \matr{F}_{\mathcal{L}_3}[3],
\end{align}
while the off-diagonal elements are
\begin{align}
&\left(\matr{M}^{-1}\right)_{21} =-\matr{F}_{{\mathcal{L}_3}\backslash\{1\}}[2]\,\matr{M}_{21}\,\matr{F}_{\mathcal{L}_3}[1],\qquad\left(\matr{M}^{-1}\right)_{12} =-\matr{F}_{{\mathcal{L}_3}\backslash\{2\}}[1]\,\matr{M}_{12}\,\matr{F}_{\mathcal{L}_3}[2],\\
&\left(\matr{M}^{-1}\right)_{32} =-\matr{F}_{{\mathcal{L}_3}\backslash\{2\}}[3]\,\matr{M}_{32}\,\matr{F}_{\mathcal{L}_3}[2],\qquad\left(\matr{M}^{-1}\right)_{23} =- \matr{F}_{{\mathcal{L}_3}\backslash\{3\}}[2]\,\matr{M}_{23}\,\matr{F}_{\mathcal{L}_3}[3],\\
&\left(\matr{M}^{-1}\right)_{31} =\matr{F}_{{\mathcal{L}_3}\backslash\{1,2\}}[3]\,\matr{M}_{32}\, \matr{F}_{{\mathcal{L}_3}\backslash\{1\}}[2]\,\matr{M}_{21}\,\matr{F}_{\mathcal{L}_3}[1],\\
&\left(\matr{M}^{-1}\right)_{13} =\matr{F}_{{\mathcal{L}_3}\backslash\{3,2\}}[1]\,\matr{M}_{12}\, \matr{F}_{{\mathcal{L}_3}\backslash\{3\}}[2]\,\matr{M}_{23}\,\matr{F}_{\mathcal{L}_3}[3].
\end{align}
\end{subequations}
The matrices $\matr{F}_\mathcal{G}[\alpha]$ are evaluated according to the recursive definition in Eq.~\eqref{eqn:MatrixInverseTheorem2}; for example
\begin{align}
  \matr{F}_{\mathcal{L}_3}[1] \hspace{-.4mm}=\hspace{-.4mm} \left[\matr{M}_1 - \matr{M}_{12}\matr{F}_{\mathcal{L}_3\backslash\{1\}}[2]\matr{M}_{21}\right]^{-1}\hspace{-1.5mm}&=\hspace{-.4mm} \left[ \matr{M}_1 - \matr{M}_{12}[\matr{M}_2-\matr{M}_{23}\matr{F}_{\mathcal{L}_3\backslash\{1,2\}}[3]\matr{M}_{32}]^{-1}\matr{M}_{21}\right]^{-1}\hspace{-1mm},\nonumber\\
  &=\hspace{-.4mm} \left[ \matr{M}_1 \hspace{-.25mm}-\hspace{-.25mm} \matr{M}_{12}[\matr{M}_2\hspace{-.25mm}-\hspace{-.25mm}\matr{M}_{23}\matr{M}_3^{-1}\matr{M}_{32}]^{-1}\matr{M}_{21}\right]^{-1}\hspace{-3mm}.\hspace{-.2mm}
\end{align}
Evaluating the flips and statics and reassembling them into matrix form gives
\begin{align}
  \matr{M}^{-1} =\frac{1}{8}\begin{pmatrix} 8 & 0 & -4 & 0 & 0 \\
 64 & -32 & -16 & -32 & 40 \\
 -16 & 16 & 4 & 16 & -16 \\
 -60 & 16 & 15 & 12 & -20 \\
 -32 & 0 & 8 & 0 & 0
  \end{pmatrix},
\end{align}
which is readily verified to be the inverse of $\matr{M}$ and is identical to the evaluation of $\matr{M}^{-1}$ as found in example \ref{exmp:GaloisMq} using Theorem \ref{thm:PathSumMatrixPower}.

\begin{exmp}[Matrix exponential]\end{exmp}As an example of the application of Theorem \ref{thm:PathSumMatrixExp} and Lemma~\ref{thm3:ResummedMatrixExponential} we consider the matrix exponential of
\begin{align}
  \matr{M}=\begin{pmatrix}
  1-\rmi&0&-\rmi&0\\0&2-\rmi&-1/3& 0\\\rmi&0&-\rmi&0\\3&-7/2&1&-1
\end{pmatrix}.
\end{align}
We use a tensor product partition 
\begin{align}
\matr{M}_{1}=\begin{pmatrix}1-\rmi &0\\0&2-\rmi\end{pmatrix},\hspace{2mm}
\matr{M}_{2}=\begin{pmatrix}-\rmi &0\\1&-1\end{pmatrix},\hspace{2mm}
\matr{M}_{12}=\begin{pmatrix}-\rmi &0\\-1/3&0\end{pmatrix},\hspace{2mm}
\matr{M}_{21}=\begin{pmatrix}\rmi &0\\3&-7/2\end{pmatrix}.\nonumber
\end{align}
Since every element of this matrix partition is non-zero, the corresponding graph is $\mathcal{K}_2$. Let us focus on the element $\exp(\matr{M})_{11}$, which forms the top-left corner of the full matrix $\exp(\matr{M})$. By Theorem \ref{thm:PathSumMatrixExp} the exact result for $\exp(t\matr{M})_{11}$ is given by
\begin{subequations}
\begin{align}
  \mathcal{L}\big[\big(\mathrm{e}^{t\matr{M}}\big)_{11}\big] = \matr{F}_{\mathcal{K}_2}[1] &= \big[s\matr{I}-\matr{M_1}-\matr{M}_{12}[s\matr{I}-\matr{M}_2]^{-1}\matr{M}_{21}\big]^{-1},\\&=\begin{pmatrix}
    \dfrac{s+\rmi}{s^2-(1-2\rmi)s-(2+\rmi)} & 0\\[4mm] \dfrac{-\rmi}{3s^3-(9-9\rmi)s^2-(6+18\rmi)s+15}&
    \dfrac{1}{s-(2 -\rmi)}
  \end{pmatrix}.
\end{align}
\end{subequations}
The inverse Laplace transform can be carried out analytically; setting $t=1$ in the result gives
\begin{subequations}
\begin{align}
  \left(\mathrm{e}^\matr{M}\right)_{11}&=\frac{\mathrm{e}^{(1-2\rmi)/2}}{15}\begin{pmatrix}
    3\left(5\cosh \tfrac{\sqrt{5}}{2}+\sqrt{5}\sinh\tfrac{\sqrt{5}}{2}\right) & 0 \\ -\rmi\left(5\mathrm{e}^{3/2}-5\cosh \tfrac{\sqrt{5}}{2}-3\sqrt{5}\sinh\tfrac{\sqrt{5}}{2}\right)& 15\,\mathrm{e}^{3/2}
  \end{pmatrix}, \\&\approx
  \begin{pmatrix}
   2.05220 - 3.19611\,\rmi& 0\\ -0.442190-0.283927\,\rmi& 3.99232-6.21768\,\rmi
  \end{pmatrix}.
\end{align}
\end{subequations}
Alternatively, we can evaluate this element by using a walk-sum, as presented in Lemma~\ref{thm3:ResummedMatrixExponential}:
\begin{align}
  \big(\mathrm{e}^\matr{M}\big)_{11}=\mathrm{e}^{\matr{M}_1} &+ \int_0^1\int_0^{t_2}\mathrm{e}^{(1-t_2)\matr{M}_1}\,\matr{M}_{12}\, \mathrm{e}^{(t_2-t_1)\matr{M_2}}\,\matr{M}_{21}\,\mathrm{e}^{t_1\matr{M}_1}\,\mathrm{d}t_1\,\mathrm{d}t_2\\
  &\hspace{-21mm}+\hspace{-1mm}\int_0^1\hspace{-2mm}\cdots\hspace{-1mm}\int_0^{t_2}\hspace{-1.5mm}
  \mathrm{e}^{(1-t_4)\matr{M}_1}\matr{M}_{12}\,\mathrm{e}^{(t_4-t_3)\matr{M}_2}\,\matr{M}_{21}\,
  \mathrm{e}^{(t_3-t_2)\matr{M}_1}\,\matr{M}_{12}\,\mathrm{e}^{(t_2-t_1) \matr{M}_2}\, \matr{M}_{21}\, \mathrm{e}^{t_1\matr{M}_1}\,\mathrm{d}t_1\hspace{-.3mm}\cdots\mathrm{d}t_4+\cdots \nonumber
\end{align}
Evaluating these terms yields
\begin{align}
  \big(\mathrm{e}^\matr{M}\big)_{11}\simeq\begin{pmatrix}
     2.05083-3.19398\,\rmi& 0\\-0.441354-0.283390\,\rmi &3.99232-6.21768\,\rmi
  \end{pmatrix}.
\end{align}
Although this result has been obtained by evaluating only the first three terms of an infinite series, it is already an excellent approximation to the exact answer: the maximum absolute elementwise error is $\sim2.5\times 10^{-3}$. This rapid convergence results from the exact resummation of the terms in the original Taylor series that correspond to walks on the graph $\mathcal{K}_2$ that contain loops.

\begin{exmp}[Matrix logarithm of a defective matrix]\label{exmp:logMat}\end{exmp}We compute the principal logarithm of the matrix
\begin{align}
&\matr{M}=\begin{pmatrix}
 4 & 1 & 1 & 2 & -1 \\
 -2 & 7 & 1 & 0 & -1 \\
 0 & -1 & 5 & 2 & 1 \\
 -2 & 0 & 0 & 8 & 0 \\
 -2 & -4 & -4 & 6 & 6
 \end{pmatrix}.
\end{align}
This matrix has only one eigenvector associated to the fifth-fold degenerate eigenvalue 6. Note also that $\|\matr{M}\|_{op}\simeq12$, with $\|.\|_{op}$ the operator norm.
We choose to partition $\matr{M}$ onto vector spaces $V_1=\mathrm{span}\big(v_1,v_3,v_4\big)$ and $V_2=\mathrm{span}\big(v_2,v_5\big)$. The corresponding graph $\mathcal{G}$ is the complete directed graph on two vertices $\mathcal{K}_2$.
Following Theorem \ref{thm:PathSumMatrixLog} we have
\begin{equation}
\matr{F}_{{\mathcal{K}_2}}[i]=\big[\matr{I}-x(\matr{I}-\matr{M}_{ii})-x^2\thinspace\matr{M}_{ij}\matr{F}_{{\mathcal{K}_2\backslash i}}[j]\matr{M}_{ji}\big]^{-1}\hspace{-1.2mm},~\textrm{and}~~\matr{F}_{{\mathcal{K}_2\backslash i}}[j]=\big[\matr{I}-x(\matr{I}-\matr{M}_{jj})\big]^{-1}\hspace{-1.2mm},\\
\end{equation}
where $(i,j)=1,2$, $i\neq j$ and $\matr{I}$ is the identity matrix of appropriate dimension.
Then the matrix elements of the principal logarithm of $\matr{M}$ are given by
\begin{equation}
\big(\log \matr{M}\big)_{ii}=\int_{0}^1\hspace{-1mm}x^{-1}(\matr{I}-\matr{F}_{{\mathcal{K}_2}}[i])\thinspace \rmd x,\quad\big(\log \matr{M}\big)_{ji}=-\int_{0}^1\matr{F}_{{\mathcal{K}_2}\backslash i}[j]\matr{M}_{ji}\matr{F}_{{\mathcal{K}_2}}[i]
\thinspace \rmd x.
\end{equation}
Reassembling them in a matrix, we obtain
\begin{equation}
\log \matr{M}=\matr{I}\, \log 6+\frac{1}{324}\begin{pmatrix}
-108 & 41 & 41 & 130 & -63 \\
 -126 & 42 & 42 & 41 & -65 \\
 18 & -37 & -37 & 71 & 56 \\
 -108 & 5 & 5 & 112 & -9 \\
 -108 & -211 & -211 & 328 & -9
\end{pmatrix}.
\end{equation}

\begin{exmp}[Matrix exponential of a quaternionic matrix]\end{exmp}The graph-theoretic nature of Theorem \ref{PSresult} allows the method to be extended to matrices over non-commutative division rings $\mathbb{D}$. In particular, 
we will show that 
Theorems \ref{thm:PathSumMatrixPower}-\ref{thm:PathSumMatrixLog} hold for quaternionic matrices $\matr{M}\in\mathbb{H}^{D\times D}$ as well. For example, consider calculating $\exp(\pi\matr{M})$ with $\matr{M}$ the quaternionic matrix
\begin{equation}
\matr{M}=\begin{pmatrix}\rmi&\mathrm{j}\\\mathrm{k}&1 \end{pmatrix},
\end{equation}
where $\mathrm{i}$, $\mathrm{j}$, and $\mathrm{k}$ are the quaternions, which satisfy $\mathrm{i}^2=\mathrm{j}^{2}=\mathrm{k}^2=\rmi\,\mathrm{j}\,\mathrm{k}=-1$. Following Theorem \ref{thm:PathSumMatrixExp}, we find for example the matrix element $\matr{M}_{11}$ to be in the Laplace domain
\begin{equation}
\mathcal{L}\big[\exp(t\matr{M})_{11}\big]=(s-\mathrm{i}-\mathrm{j}(s-1)^{-1}\mathrm{k})^{-1}.
\end{equation}
Calculating the other matrix elements, inverting the Laplace transforms and setting $t=\pi$, we obtain
\begin{equation}
\exp(\pi\matr{M})=-\frac{1+e^{\pi}}{2}\begin{pmatrix} \mathrm{i}+ \tanh \left(\frac{\pi }{2}\right) & \mathrm{j}-\mathrm{k} \\
 \mathrm{k}-\mathrm{j} & \mathrm{i}+ \tanh \left(\frac{\pi }{2}\right)\end{pmatrix}.
\end{equation}

\subsection{\textbf{Exact matrix exponential of block tridiagonal matrices}}\label{subsec:ExmpLnExp}
Let $\{\matr{M}_{k',k}\}$ be a partition of a matrix $\matr{M}$ such that $\matr{M}_{k'\neq k\pm 1,k}=0$.  If this partition consists of $N^2$ pieces, the corresponding graph is the finite linear graph $\mathcal{L}_{N}$, and $\matr{M}$ can therefore be said to be an $\mathcal{L}_N$-structured matrix. 
As mentioned in remarks \ref{rem:block} and \ref{rem:KnownFormulae}, such a matrix is essentially a block tridiagonal matrix since there exists a permutation matrix $\matr{P}$ such that $\matr{P}\matr{M}\matr{P}^{\mathrm{T}}$ is block tridiagonal. 
For these matrices, the path-sum expression for the matrix exponential and inverse can be written in a particularly compact form. For $k=1,\ldots,N$, we set $\tilde{\matr{M}}_{k}=s\matr{I}-\matr{M}_k$, and define the finite continued fractions
\begin{align}
\tilde{\matr{X}}_k=\Big[\tilde{\matr{M}}_{k}-\matr{M}_{k,k+1}\tilde{\matr{X}}_{k+1}\matr{M}_{k+1,k}\Big]^{-1},\quad
\tilde{\matr{Y}}_k=\Big[\tilde{\matr{M}}_{k}-\matr{M}_{k,k-1}\tilde{\matr{Y}}_{k-1}\matr{M}_{k-1,k}\Big]^{-1},\label{eqn:Xk}
\end{align}
with $\tilde{\matr{X}}_{N}\equiv\tilde{\matr{M}}_N^{-1}$ and $\tilde{\matr{Y}}_{1}\equiv\tilde{\matr{M}}_1^{-1}$. Note that the inversion height of $\tilde{\matr{X}}_k$ $(\tilde{\matr{Y}}_k)$ is $N+1-k$ $(k)$. Let $\tilde{\matr{U}} = \mathfrak{L}[\exp(t\,\matr{M})]$ be the Laplace transform of the matrix exponential of $\matr{M}$; then the partition of $\tilde{\matr{U}}$ is given by
\begin{subequations}
\begin{gather}
\tilde{\matr{U}}_{kk}=\Big[\,\tilde{\matr{X}}^{-1}_k+\tilde{\matr{Y}}_k^{-1}-\tilde{\matr{M}}_k\,\Big]^{-1},\label{eqn:UkonLN}\\
\tilde{\matr{U}}_{k,k'<k}=\hspace{-2mm}\prod_{j=k'+1}^{k}\hspace{-1mm}\big(\tilde{\matr{X}}_j \matr{M}_{j,j-1}\big)\,\,\tilde{\matr{U}}_{k'k'},\qquad\tilde{\matr{U}}_{k,k'>k}=\hspace{-2mm}\prod_{j=k'-1}^{k}\hspace{-1mm}\big(\tilde{\matr{Y}}_j \matr{M}_{j,j+1}\big)\,\,\tilde{\matr{U}}_{k'k'},\label{eqn:UkkprimeonLN}
\end{gather}
\end{subequations}
where $\prod_{j=1}^na_j=a_na_{n-1}\cdots a_2a_1$ is a left product. Provided the required inverse Laplace-transforms are analytically available, Eq.(\ref{eqn:Xk}-\ref{eqn:UkkprimeonLN}) yields the exact matrix exponential of $\matr{M}$.
Similar formulae for the inverse of arbitrary $\mathcal{L}_N$-structured matrices are obtained upon replacing $s\matr{I}$ by $\matr{I}$ in each of the $\tilde{\matr{M}}_k$.


\subsection{\textbf{Computational cost of path-sum on arbitrary tree-structured matrices}}\label{subsec:ExmpTree}
The number of paths in any finite graph is finite. Thus, as soon as a matrix $\matr{M}$ has a finite size, the path-sum expressions of Theorems \ref{thm:PathSumMatrixPower}-\ref{thm:PathSumMatrixLog} evaluate any function $f(\matr{M})$ exactly in a finite number $\rho$ of operations. In this final example we calculate $\rho$ in the case of matrices $\matr{M}$ whose partition is an arbitrary tree, denoted by $\mathcal{T}_N$.

First, we consider the computational cost of evaluating a static $f(\matr{M})_{\alpha\alpha}$. Since a tree contains no bare cycles of length greater than 2, the vertex $\alpha$ requires dressing only by loops and edge cycles on $\mathcal{T}_N$. Consequently, the sequence of operations involved in dressing $\alpha$ fall into two categories: nesting (adding an edge cycle which dresses the internal vertex of a previous edge cycle) or branching (including an extra edge cycle at constant dressing depth). These two operations have the same computational cost, as each requires one inversion, two multiplications and one subtraction of $d\times d$ matrices. Due to this symmetry between branching and nesting, the computational costs of dressing a vertex on the linear graph on $N$ vertices $\mathcal{L}_N$ (which involves nesting only) and the star graph on $N$ vertices $\mathcal{S}_N$ (which involves branching only) are identical: each requires $N$ inversions, $2(N-1)$ multiplications and $N-1$ additions to fully dress any vertex.
Since any tree can be decomposed as an ensemble of linear graphs (the branches) and star graphs (the nodes), the computational cost of evaluating $f(\matr{M})_{\alpha\alpha}$ when $\matr{M}$ is tree-structured depends only on the number of vertices of the tree: the cost is independent of the detailed structure of $\mathcal{T}_N$. The number of floating point operations required to evaluate the diagonal element $f(\matr{M})_{\alpha\alpha}$, denoted by $\rho_{\mathcal{T}_N;\alpha}$, scales as $3Nd^3$, while the computational density (the number of operations per usual matrix-element evaluated) is $\sim3Nd$.

We now consider the cost of calculating an off-diagonal element $f(\matr{M})_{\omega\alpha}$. Let $(\alpha\nu_2\cdots\nu_{\ell}\omega)$ be the unique path leading from $\alpha$ to $\omega$ on $\mathcal{T}_N$, and $\rho_{\mathcal{T}_N;\alpha}$ be the cost of fully dressing the vertex $\alpha$ on $\mathcal{T}_N$. Then the total cost of evaluating $f(\matr{M})_{\omega\alpha}$ is
\begin{align}\label{eqn:TreeCompCost}
2\ell d^3 + \sum_{i=1}^{\ell+1} \rho_{\mathcal{T}_N\backslash\{\alpha,\ldots,\nu_{i-1}\};\nu_i},
\end{align}
where the first term accounts for the cost of evaluating the matrix multiplications along the path. Note that this cost depends on the structure of $\mathcal{T}_N$, since the number of vertices in each tree in the sequence $\mathcal{T}_N\backslash\{\alpha,\ldots,\nu_{i-1}\}$ depends on whether any of the previously-visited vertices are nodes of $\mathcal{T}_N$. Nevertheless, we can place an upper bound on the cost by considering the case where the tree $\mathcal{T}_N$ is the linear graph $\mathcal{L}_N$. Then $\mathcal{L}_N\backslash\{\alpha,\ldots,\nu_{i-1}\}$ contains $N+1-i$ vertices, and the total cost of evaluating the contribution of a path of length $\ell$ is  $\sim \ell d^3(3N+2)$.
Finally, we note that in the course of evaluating the contribution of a certain path $p$, we simultaneously evaluate the contributions of all subpaths of $p$. The operations counted by Eq.~\eqref{eqn:TreeCompCost} therefore generate $(\ell+1)d^2$ elements of the partition of $f\big(\matr{M}\big)$. This improves the computational density of the path-sum method for tree-structured matrices to $\sim 3Nd$. Therefore, we conclude that evaluating any element of a partition of a matrix function of a tree-structured matrix using path-sums is efficient: i.e.~the cost is linear in the number of vertices in the tree.

\begin{figure}[t!]
\begin{center}
\vspace{-2mm}
\includegraphics[width=.5\textwidth]{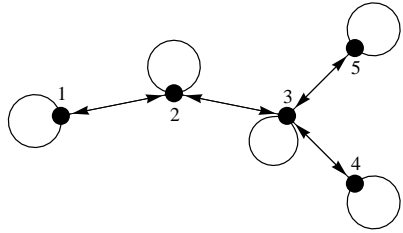}
\vspace{-5mm}
\caption{A tree on $N=5$ vertices. All edges are bidirectional.}
\label{fig:TreeCost}
\end{center}
\vspace{-9.5mm}
\end{figure}

\begin{exmp}[Vertex dressing on a tree]\end{exmp}To illustrate the above discussion, consider the cost of calculating the inverse of a matrix $\matr{M}\in\mathbb{C}^{5d\times 5d}$, $d\in\mathbb{N}^*$ whose tensor-product partition into $d\times d$ blocks yields the tree on $N=5$ vertices represented on Fig. \ref{fig:TreeCost}. Let $\matr{M}_{\mu\nu}$ be the submatrix associated to the edge $(\nu\mu)$ on the graph.
Now consider calculating $(\matr{M}^{-1})_{3}$. According to Theorem \ref{thm:PathSumMatrixInverse}, it is
\begin{align}
\hspace{-0mm}(\matr{M}^{-1})_{3}=\Big[\hspace{-.5mm}\matr{M}_3&-\overbrace{\matr{M}_{32}[\underbrace{\matr{M}_2-\matr{M}_{21}\matr{M}_1^{-1}\matr{M}_{12}}_{\textrm{nesting}}]^{-1}\matr{M}_{23}}^{\textrm{branching}}-\overbrace{\matr{M}_{34}\matr{M}_4^{-1}\matr{M}_{43}}^{\textrm{branching}}-\overbrace{\matr{M}_{35}\matr{M}_5^{-1}\matr{M}_{53}}^{\textrm{branching}}\hspace{-.2mm}\Big]^{-1}.\nonumber
\end{align}
This comprises three branching operations corresponding to the three branches attached to vertex (3) and one nesting corresponding to the dressing of vertex (2) by cycles through (1). In total $(\matr{M}^{-1})_{3}$ thus requires 5 inversions, 8 multiplications and 3 additions of $d\times d$ matrices, that is $(5+8)d^3+3d^2$ operations which scales as $3d^3 N$ as expected. Now consider calculating $(\matr{M}^{-1})_{1}$. Following Theorem \ref{thm:PathSumMatrixInverse} we find
\begin{align}
\hspace{-1mm}(\matr{M}^{-1})_{1}=\hspace{-.5mm}\Big[\hspace{-.5mm}\matr{M}_1&-\matr{M}_{12}\big[\matr{M}_2-\matr{M}_{23}[\matr{M}_3-\matr{M}_{34}\matr{M}_4^{-1}\matr{M}_{43}-\matr{M}_{35}\matr{M}_5^{-1}\matr{M}_{53}]^{-1}\matr{M}_{32}\big]^{-1}\matr{M}_{21}\hspace{-.5mm}\Big]^{-1}.\nonumber
\end{align}
This comprises two successive nesting operations corresponding to the dressing of vertex (1) by cycles through (2) and (2) by cycles through (3) and two branching operations corresponding to remaining branches off (3). In total $(\matr{M}^{-1})_{1}$ thus requires 5 inversions, 8 multiplications and 3 additions of $d\times d$ matrices exactly like $(\matr{M}^{-1})_{3}$, as expected. The computational cost of calculating any element of $\matr{M}^q$, $\exp(\matr{M})$ or $\log(\matr{M})$ is identical, except for the additional costs associated to the inverse Z-transform, Laplace-transform and Richter integral, respectively.

\section{Conclusion and Outlook}\label{sec:V:Conclusion}
The method of path-sums is based on three main concepts: firstly, the partitioning of a discrete matrix $\matr{M}$ into an ensemble of submatrices whose dimensionalities can be freely chosen; secondly, the mapping between multiplication of $\matr{M}$ and walks on a graph $\mathcal{G}$ whose edge pattern encodes the structure of $\matr{M}$; and thirdly, the exact closed-form resummation of certain classes of walks on $\mathcal{G}$ through the dressing of vertices by cycles. By combining these concepts, any partition of a function of a finite matrix $f\big(\matr{M}\big)$ 
can be exactly evaluated in a finite number of steps, provided the required inverse transforms are analytically available.
The graph-theoretic nature of Theorem \ref{PSresult} permits the extension of the method of path-sums to functions of continuous operators as well as matrices with entries in non-commutative division rings; as we will show elsewhere. 
%

Using a directed graph to encode the structure of a matrix allows any structure and symmetries that the matrix possesses to be easily recognized and exploited. We thus expect the method of path-sums to have widespread applications in, for example, the study of Markov chains and quantum many-body physics, where the relevant matrix -- the many-body Hamiltonian -- is both sparse and highly structured. We have already successfully implemented the method to study the dynamics of Rydberg-excited Mott insulators \cite{Giscard2011a}. 
We also note that a generalization of the walk-sum result for the matrix exponential to continuous matrices is already known, and leads to Feynman path integrals, while the cycle dressing of a vertex yields the Dyson equation; these results will be presented elsewhere. We expect that this approach to path integrals will help understand their divergences.

Finally, the development and study of numerical implementations of the path-sum method is a topic of key importance for future work.  
Indeed for most applications requiring the computation of a function of a matrix, an inexact but accurate result is sufficient, in particular if the matrix $\matr{M}$ itself is known only with a finite precision. Additionally on highly connected graphs the number of paths and bare cycles can be very large, so that an exact evaluation of the flips and statics of $f\big(\matr{M}\big)$ may have a prohibitive computational cost. In these cases, it is desirable to truncate the exact expressions given in Theorems \ref{thm:PathSumMatrixPower}-\ref{thm:PathSumMatrixLog}, e.g. by dressing vertices to a depth smaller than their maximum dressing depth or by neglecting certain paths. This may also be motivated by considerations external to the method: in quantum mechanics, paths and cycles on the graph of the Hamiltonian represent physical processes (similarly to Feynman diagrams) some of which might be negligible. The accuracy of truncated numerical approximations to the path-sum will thus be of paramount importance in many applications.
Recent research \cite{BenziNEW, Benzi1999, Benzi2007, Cramer2006, Iserles2000} has shown that the entries $f\big(\matr{M}\big)_{ij}$ of certain matrix functions decay exponentially with the length of the shortest path between the vertices that correspond to entries $\matr{M}_{ii}$ and $\matr{M}_{jj}$. If these results can be extended to the norm of the statics and flips of a general partition of $f\big(\matr{M}\big)$, they could be used to estimate the contribution of each path in a path-sum expression, thereby giving some indication of how the truncation of a path-sum expression would affect the accuracy of the result.

\section*{Acknowledgments}
We thank M. Benzi for fruitful discussions.

\Appendix
\section{Proofs of Path-Sum Expressions for Common Matrix Functions}\label{appendix:MatrixFunctionsProofs}
In this appendix we prove the path-sum results presented earlier without proof. We begin by proving the results of \S \ref{subsec:PartitionOfMatrixPowers} relating the partitions of $\matr{M}^k$ and $f(\matr{M})$ to the partition of $\matr{M}$. We then prove Theorems \ref{thm:PathSumMatrixPower}-\ref{thm:PathSumMatrixLog}; i.e.~the path-sum results for a matrix raised to a complex power (\S\ref{subsec:MatrixPowerProof}), the matrix inverse (\S\ref{subsec:MatrixInverseProof}), the matrix exponential (\S\ref{subsec:MatrixExpProof}), and the matrix logarithm (\S\ref{subsec:MatrixLogProof}). 

\subsection{\textbf{Partitions of matrix powers and functions}}\label{subsec:MatrixFunctionProof}
Consider an element $\matr{R}_{\omega}\matr{M}^k\matr{R}^{\mathrm{T}}_{\alpha}$ of a general partition of $\matr{M}^k$. On inserting the identity in the form of the closure relation over projector-lattices between each appearance of $\matr{M}$ in the product and using $\varepsilon_\mu=\matr{R}_\mu^{\mathrm{T}}\matr{R}_\mu$, we obtain
\begin{align}
  \left(\matr{M}^k\right)_{\omega\alpha} =\hspace{-1.5mm} \sum_{\eta_{k},\ldots,\eta_2}\matr{R}_{\omega}\matr{M}\, \matr{R}^{\mathrm{T}}_{\eta_{k}}\matr{R}_{\eta_{k}}\cdots \matr{M}\, \matr{R}^{\mathrm{T}}_{\eta_{2}}\matr{R}_{\eta_{2}} \matr{M}\,\matr{R}^{\mathrm{T}}_\alpha=\hspace{-1.5mm}\sum_{\eta_{k},\ldots,\eta_2}\matr{M}_{\omega\eta_{k}} \cdots \matr{M}_{\eta_3\eta_2} \matr{M}_{\eta_2\alpha},
\end{align}
with $\alpha\equiv\eta_1$ and $\omega\equiv \eta_{k+1}$. This expression describes matrix multiplication in terms of the partition of $\matr{M}$ and provides an explicit description of which pieces of $\matr{M}$ contribute to a given piece of $\matr{M}^k$. It follows that the partition of a matrix function $f(\matr{M})$ with power series $f\big(\matr{M}\big) = \sum_{k=0}^\infty f_k\, \matr{M}^k$ is given by
\begin{align}\label{eqn:MatrixFunctionDecompApp}
  f\big(\matr{M}\big)_{\omega\alpha} = \sum_{k=0}^\infty f_k\,\sum_{\eta_{k},\ldots,\eta_2}\matr{M}_{\omega\eta_{k}} \cdots \matr{M}_{\eta_3\eta_2} \matr{M}_{\eta_2\alpha}.
\end{align}
This equation relates the partition of $f\big(\matr{M}\big)$ to that of $\matr{M}$.

\subsection{\label{subsec:MatrixPowerProof}\textbf{A matrix raised to a complex power}}
\begin{proof} We consider a matrix $\matr{M}\in\mathbb{C}^{D\times D}$.
To prove Theorem \ref{thm:PathSumMatrixPower} we start from the power series
\begin{align}\label{eqn:ComplexPowerAsPowerSeries}
\matr{M}^q = \sum_{k=0}^\infty \binom{q}{k}(\matr{M}-\matr{I})^k,
\end{align}
where $q\in\mathbb{C}$ and $\tbinom{q}{k}=q^{\underline{k}}/k!$ is a binomial coefficient, with $q^{\underline{k}}$ the falling factorial. Note that although the sum only converges when $\|\matr{M}-\matr{I}\| <1$, with $\|.\|$ a sub-multiplicative norm, the result of Theorem \ref{thm:PathSumMatrixPower} is valid for matrices of arbitrary norm by analytic continuation. By applying the result of Eq.~(\ref{eqn:MatrixFunctionDecomp}) to Eq.~(\ref{eqn:ComplexPowerAsPowerSeries}) we find that an element of a partition of $\matr{M}^q$ is given by
\begin{equation}
\big(\matr{M}^q\big)_{\omega\alpha}=\sum_{k=0}^{\infty}\binom{q}{k}\sum_{W_{\mathcal{G};\alpha\omega;k}}\bar{\matr{M}}_{\omega\eta_k}\cdots\bar{\matr{M}}_{\eta_3\eta_2}\bar{\matr{M}}_{\eta_2\alpha},
\end{equation}
where we have introduced the auxiliary matrix $\bar{\matr{M}} = \matr{M}-\matr{I}$ and $\mathcal{G}$ is the graph of the partition of $\bar{\matr{M}}$. We shall now recast this expression so as to make the loops of the walk $(\alpha)(\alpha\eta_2)\cdots(\eta_k\omega)(\omega)$ appear explicitly. To this end, we remark that when a loop off a vertex $\mu$ occurs $p\in\mathbb{N}$ consecutive times in a walk, the contribution of the walk is comprises a factor of $(\bar{\matr{M}}_\mu)^p$. Thus
\begin{align}
&\hspace{-2.5mm}\big(\matr{M}^q\big)_{\omega\alpha}\hspace{-.2mm}=\hspace{-.1mm}\sum_{k=0}^{\infty}\hspace{-.5mm}\binom{q}{k}\hspace{-.5mm}\sum_{m=0}^k\sum_{W_{\mathcal{G}_0;\alpha\omega;m}}\sum_{\{p_i\}/ k-m}\hspace{-1.5mm}(\bar{\matr{M}}_{\omega})^{p_{m+1}}\matr{M}_{\omega \mu_m}\hspace{-.3mm}\cdots(\bar{\matr{M}}_{\mu_2})^{p_2}\matr{M}_{\mu_2\alpha}(\bar{\matr{M}}_{\alpha})^{p_1}.
\end{align}
In this expression any two consecutive vertices $\mu_\ell$ and $\mu_{\ell+1}$ are now distinct and the sum $\sum_{W_{\mathcal{G}_0;\alpha\omega;m}}$ runs thus over the walks of the loopless graph $\mathcal{G}_0$. Each integer $p_i\in\mathbb{N}$, called a loop number, represents the number of consecutive times a loop is undergone off a certain vertex. The notation $\{p_i\}/ k-m$ on the final sum therefore indicates that the final sum runs over all possible configurations of the $p_i$ numbers, that is all configurations of $k-m$ loops on $m+1$ vertices, subject to the restriction that the loop number $p_i$ on any loopless vertex is fixed to zero. This implies that $\sum_{i=1}^{m+1}p_i=k-m$ and the $p_i$ are thus said to form a weak composition of $k-m$.
We now remark that for any such weak composition and $q\in\mathbb{C}$, the following relation holds:
\begin{align}
\label{eq:BinomialComposition}
&\binom{q}{k}=\sum_{k_m}\sum_{k_{m-1}=0}^{k_m}\cdots\sum_{k_1=0}^{k_2}\frac{(-1)^m}{\prod_{r=1}^{m+1}p_r!}(-1)^{p_{m+1}}\big((-q-1)-k_m+p_{m+1}\big)^{\underline{p_{m+1}}}\cdots\nonumber\\
&\hspace{4.8cm}\cdots(-1)^{p_2}(k_2-k_1+p_2)^{\underline{p_2}}\,\,(-1)^{p_1}(k_1+p_1)^{\underline{p_1}}\,,
\end{align}
where the first sum is an \emph{indefinite} sum to be evaluated at $k_m=-q-1$. From now on we will denote this by $\sum_{k_m=0}^{-q-1}$. Equation \eqref{eq:BinomialComposition}, which is independent of the value of each individual $p_i$, is proved by induction on $m$. Upon substituting this relation into Eq.~(\ref{eqn:ComplexPowerAsPowerSeries}) and rearranging the order of the summations we obtain
\begin{align}
&\big(\matr{M}^q\big)_{\omega\alpha}\hspace{-.6mm}=\hspace{-1.1mm}\sum_{m=0}^\infty(-1)^m\hspace{-4mm}\sum_{W_{\mathcal{G}_0;\alpha\omega;m}}\hspace{-.5mm}\sum_{\{p_i\}=0}^{\infty}\sum_{k_m=0}^{-q-1}\sum_{k_{m-1}=0}^{k_m}\hspace{-1mm}\cdots\hspace{-1mm}\sum_{k_1=0}^{k_2}\hspace{-.5mm}(-\bar{\matr{M}}_{\omega})^{p_{m+1}}\big((-q\hspace{-.2mm}-\hspace{-.5mm}1)\hspace{-.5mm}-\hspace{-.5mm}k_m\hspace{-.5mm}+\hspace{-.5mm}p_{m+1}\big)^{\underline{p_{m+1}}}\matr{M}_{\omega \mu_m}\cdots\nonumber\\
&\hspace{3cm}\cdots(-\bar{\matr{M}}_{\mu_2})^{p_{2}}(k_2-k_1+p_2)^{\underline{p_2}}\matr{M}_{\mu_2\alpha}(\bar{\matr{M}}_{\alpha})^{p_1}(k_1+p_1)^{\underline{p_1}}\,,\label{eq:MatrixPowerNoResum}
\end{align}
where the sum $\sum_{\{p_i\}=0}^{\infty}=\sum_{p_1=0}^\infty\cdots \sum_{p_{m+1}=0}^\infty$ runs over all the loop numbers, subject to the restriction that the loop number $p_i$ on any loopless vertex is fixed to zero. We now evaluate the contributions from the infinite loop sums in closed form by noting that
\begin{equation}
\label{eq:LoopSumMatrixPower}
\sum_{p_i=0}^\infty (-\bar{\matr{M}}_{\mu})^{p_i}\frac{(k_i-k_{i-1}+p_i)^{\underline{p_i}}}{p_i!}=(\matr{M}_{\mu})^{-(k_i-k_{i-1})-1},
\end{equation}
where we have used the fact that $\|\bar{\matr{M}}_{\mu}\|= \|\matr{M}_{\mu}-\matr{I}\| \le \|\matr{M}-\matr{I}\|<1$. Introducing Eq.~(\ref{eq:LoopSumMatrixPower}) into Eq.~(\ref{eq:MatrixPowerNoResum}) and setting $k_0 = 0$ by definition yields the expression
\begin{align}
&\big(\matr{M}^q\big)_{\omega\alpha}=\sum_{m=0}^\infty(-1)^m\sum_{W_{\mathcal{G}_0;\alpha\omega;m}}\sum_{k_m=0}^{-q-1}\sum_{k_{m-1}=0}^{k_m}\cdots\sum_{k_1=0}^{k_2}(\matr{M}_{\omega})^{-((-q-1)-k_m)-1}\matr{M}_{\omega \mu_m}\cdots\nonumber\\
&\hspace{5.4cm}\cdots(\matr{M}_{\mu_2})^{-(k_2-k_1)-1}\matr{M}_{\mu_2\alpha}(\matr{M}_{\alpha})^{-k_1-1}.
\end{align}
This expression is an $m$-fold nested discrete convolution. In order to convert the discrete convolution to a product, we take the unilateral Z-transform of the above expression with respect to $n\equiv-q-1$. We obtain
\begin{align}
&\big(\matr{M}^q\big)_{\omega\alpha}=-\mathcal{Z}^{-1}\left\{\sum_{m=0}^\infty\sum_{W_{\mathcal{G}_0;\alpha\omega;m}}\big[\matr{I}z^{-1}-\matr{M}_{\omega}\big]^{-1}\matr{M}_{\omega \mu_m}\cdots\right.\nonumber\\
&\hspace{3cm}\left.\cdots\big[\matr{I}z^{-1}-\matr{M}_{\mu_2}\big]^{-1}\matr{M}_{\mu_2\alpha}\big[\matr{I}z^{-1}-\matr{M}_{\alpha}\big]^{-1}\vphantom{\sum_{m=0}^\infty}\right\}[n]\Big|_{n=-q-1}, \label{eqn:ComplexPowerWalkSum}
\end{align}
where $z \in \mathbb{C}$ is the Z-domain variable. Now the content of the inverse Z-transform is a sum of walk  contributions. Indeed we can see $c^{\textrm{eff}}_\alpha=\big[\matr{I}z^{-1}-\matr{M}_{\alpha}\big]^{-1}$ as an effective weight associated to vertex $\alpha$. This effective weight results from the dressing of $\alpha$ by all the loops off $\alpha$ which is performed by Eq.~(\ref{eq:LoopSumMatrixPower}). Then upon remarking that e.g. $\matr{M}_{\mu_2\alpha}=c_{\mu_2\alpha}$ is the weight associated to edge $(\alpha\mu_2)$ from $\alpha$ to $\mu_2$, Eq.~(\ref{eqn:ComplexPowerWalkSum}) is
\begin{align}
\big(\matr{M}^q\big)_{\omega\alpha}=-\mathcal{Z}^{-1}\Big\{\sum_{W_{\mathcal{G}_0;\alpha\omega}}\hspace{-1mm}c^{\textrm{eff}}_\omega\, c_{\omega\mu_m}\cdots c^{\textrm{eff}}_{\mu_2}\,c_{\mu_2\alpha}\,c^{\textrm{eff}}_{\alpha}
\Big\}[n]\Big|_{n=-q-1}. 
\end{align}
This is now in a form suitable to the use of Theorem \ref{PSresult}. We obtain
\begin{equation}
\hspace{-0mm}(\matr{M}^{q})_{\omega\alpha}=-\mathcal{Z}^{-1}\hspace{-1mm}\left\{\sum_{P_{\mathcal{G}_0;\alpha\omega}}\matr{F}_{\mathcal{G}\backslash\{\alpha,\ldots,\nu_\ell\}}[\omega]\,{\matr{M}}_{\omega\nu_\ell}\,\cdots \,\matr{F}_{\mathcal{G}\backslash\{\alpha\}}[\nu_2]\,{\matr{M}}_{\nu_2\alpha}\,\matr{F}_\mathcal{G}[\alpha]\right\}\hspace{-.5mm}[n]\Bigg|_{n=-q-1}\hspace{-8mm},\label{PathMqProof}\hspace{-10mm}
\end{equation}
where $\mathcal{G}$ is the graph of $\{(\matr{M}-\matr{I})_{\mu\nu}\}$, $\ell$ is the length of the path, and
\begin{equation}
\matr{F}_{\mathcal{G}}[\alpha] = \Big[\matr{I}-\hspace{-1mm}\sum_{C_{\mathcal{G};\alpha}}\matr{M}_{\alpha\mu_m} \matr{F}_{\mathcal{G}\backslash{\{\alpha,\ldots,\mu_{m-1}}\}}[\mu_m]\cdots \matr{F}_{\mathcal{G}\backslash\{\alpha\}}[\mu_2]\matr{M}_{\mu_2\alpha}\thinspace\Big]^{-1}\hspace{-2mm},\label{eq:Fdemo1}
\end{equation}
with $m$ the length of the bare cycle. The quantity $\matr{F}_{\mathcal{G}}[\alpha]$ can itself be seen as an effective weight associated to vertex $\alpha$ resulting from the dressing of $\alpha$ by all the closed walks off $\alpha$ in $\mathcal{G}$. Since a loop is a bare cycle, the dressing of the vertices by their loops is included in $\matr{F}_{\mathcal{G}}[\alpha]$ as well. This is obvious if one considers a
graph $\mathcal{G}$ that is reduced to a unique vertex $\alpha$ presenting a loop. In that case Eq.~(\ref{eq:Fdemo1}) yields $\matr{F}_\alpha[\alpha]=[\matr{I}z^{-1}-\matr{M}_\alpha]^{-1}\equiv c^{\textrm{eff}}_\alpha$. For convenience we can make the loop dressing completely explicit in Eq.~(\ref{eq:Fdemo1}) by separating the loops from the other bare cycles
\begin{equation}
\hspace{0mm}\matr{F}_{\mathcal{G}}[\alpha] = \Big[\matr{I}z^{-1} -\matr{M}_\alpha-\hspace{-1mm}\sum_{C_{\mathcal{G}_0;\alpha}}\matr{M}_{\alpha\mu_m} \matr{F}_{\mathcal{G}\backslash{\{\alpha,\ldots,\mu_{m-1}}\}}[\mu_m]\cdots \matr{F}_{\mathcal{G}\backslash\{\alpha\}}[\mu_2]\matr{M}_{\mu_2\alpha}\thinspace\Big]^{-1},\label{eq:Fdemo2}
\end{equation}
and note that now the sum runs over the bare cycles of the loopless graph $\mathcal{G}_0$. Together with Eq.~(\ref{PathMqProof}) the above Eq.~(\ref{eq:Fdemo2}) proves Theorem \ref{thm:PathSumMatrixPower}.
\qquad\end{proof}

\subsection{\textbf{The matrix inverse}}\label{subsec:MatrixInverseProof}
\begin{proof} We consider an invertible matrix $\matr{M}\in\mathbb{C}^{D\times D}$.
To prove Theorem \ref{thm:PathSumMatrixInverse} we write the matrix inverse as $\matr{M}^{-1} = \sum_{n=0}^\infty\left(\matr{I}-\matr{M}\right)^n$. Note that the sum only converges for $\|\matr{I}-\matr{M} \|<1$; nevertheless, the end result can be extended to matrices of arbitrary norm by analytic continuation, and Theorem \ref{thm:PathSumMatrixInverse} is therefore valid for all matrices, regardless of norm. Introducing an auxiliary matrix $\bar{\matr{M}} \equiv \matr{I}-\matr{M}$, we apply the result of Eq.~\eqref{eqn:MatrixFunctionDecomp} to the power series to obtain
\begin{align}
  \left(\matr{M}^{-1}\right)_{\omega\alpha} &= \sum_{n=0}^\infty \sum_{W_{\mathcal{G};\alpha\omega;n}} \bar{\matr{M}}_{\omega\eta_n}\cdots\bar{\matr{M}}_{\eta_3\eta_2}\bar{\matr{M}}_{\eta_2\alpha}.
\end{align}
where $\mathcal{G}$ is the graph of the partition of $\bar{\matr{M}}$. We follow the same procedure as in \S\ref{subsec:MatrixPowerProof} above and we omit the details; the result is
\begin{subequations}
\begin{align}\label{eqn:MatrixInverseWalkForn}
 \hspace{-2mm} \left(\matr{M}^{-1}\right)_{\omega\alpha}&= \sum_{n=0}^\infty\sum_{W_{\mathcal{G}_0;\alpha\omega;n}} \left[\matr{I}-\bar{\matr{M}}_{\omega}\right]^{-1}\, \bar{\matr{M}}_{\omega\nu_n} \cdots \left[\matr{I}-\bar{\matr{M}}_{\nu_2}\right]^{-1}\,\bar{\matr{M}}_{\nu_2\alpha}\,\left[\matr{I}-\bar{\matr{M}}_{\alpha}\right]^{-1},\\
 &= \sum_{n=0}^\infty\sum_{W_{\mathcal{G}_0;\alpha\omega;n}}\hspace{-1mm}(-1)^n\, \matr{M}_{\omega}^{-1}\, \matr{M}_{\omega\nu_n} \cdots \matr{M}_{\nu_2}^{-1}\,\matr{M}_{\nu_2\alpha}\,\matr{M}_{\alpha}^{-1},\label{eqn:MatrixInverseWalkSum}
\end{align}
\end{subequations}
where we have used $\bar{\matr{M}}_\mu=\matr{I}-\matr{M}_\mu$ and $\bar{\matr{M}}_{\mu \nu}=-\matr{M}_{\mu\nu}$.
Eq.~(\ref{eqn:MatrixInverseWalkSum}) is a sum of walk contributions with effective vertex weights $c^{\textrm{eff}}_\mu=\matr{M}_\mu^{-1}$ resulting from the loop dressing of $\mu$ which occurs when $\matr{M}_\mu\neq0$. We now use Theorem \ref{PSresult} and obtain
\begin{subequations}
\begin{align}
&\big(\matr{M}^{-1}\big)_{\omega\alpha}=\sum_{P_{\mathcal{G}_0;\alpha\omega}} (-1)^\ell\, \matr{F}_{\mathcal{G}\backslash\{\alpha,\ldots,\nu_\ell\}}[\omega]\matr{M}_{\omega\nu_\ell}\cdots \matr{F}_{\mathcal{G}\backslash\{\alpha\}}[\nu_2]\,\matr{M}_{\nu_2\alpha}\,\matr{F}_{\mathcal{G}}[\alpha],\\
\shortintertext{where $\mathcal{G}$ is the graph of $\{(\matr{M}-\matr{I})_{\mu\nu}\}$, $\ell$ is the length of the path, and}
&\hspace{-1mm} \matr{F}_{\mathcal{G}}[\alpha] = \Big[\matr{M}_{\alpha} - \sum_{C_{\mathcal{G}_0;\alpha}}(-1)^m\,\matr{M}_{\alpha\mu_m} \matr{F}_{\mathcal{G}\backslash{\{\alpha,\ldots,\mu_{m-1}}\}}[\mu_m]\cdots \matr{F}_{\mathcal{G}\backslash\{\alpha\}}[\mu_2]\matr{M}_{\mu_2\alpha}\thinspace\Big]^{-1}\hspace{-1.5mm},\label{eq:Finverse}
\end{align}
\end{subequations}
with $m$ the length of the bare cycle. Similarly to \S\ref{subsec:MatrixPowerProof}, we have separated the contribution of the loops from that of the other bare cycles.
This proves Theorem \ref{thm:PathSumMatrixInverse}.
\qquad\end{proof}

\subsection{\textbf{The matrix exponential}}\label{subsec:MatrixExpProof}
\begin{proof} We consider a matrix $\matr{M}\in\mathbb{C}^{D\times D}$.
To prove Theorem \ref{thm:PathSumMatrixExp} and Lemma \ref{thm3:ResummedMatrixExponential}, we start from the power series expression $\exp(t\, \matr{M}) = \sum_{n=0}^\infty t^n\matr{M}^n/n!$. By applying the result of Eq.~\eqref{eqn:MatrixFunctionDecomp} to this series we find that the partition of $\exp(t\,\matr{M})$ is given by
\begin{align}
\exp(t\,\matr{M})_{\omega\alpha} = \sum_{n=0}^\infty \frac{t^n}{n!}\sum_{W_{\mathcal{G};\alpha\omega;n}} \matr{M}_{\omega\eta_n}\cdots\matr{M}_{\eta_3\eta_2}\matr{M}_{\eta_2\alpha}.
\end{align}
with $\mathcal{G}$ the graph of the partition of $\matr{M}$.
Following the same procedure as in \S\ref{subsec:MatrixPowerProof}, we make the loop of the walks appear explicitly
\begin{multline}\label{eqn:DecompOfMatrixExpI}
 \hspace{-4mm} \exp(t\,\matr{M})_{\omega\alpha} =\hspace{-.7mm}\sum_{n=0}^{\infty} \frac{t^n}{n!}\sum_{m=0}^n \sum_{W_{\mathcal{G}_0;\alpha\omega;m}}\sum_{\{p_i\}/ n-m}\hspace{-3mm} \left(\matr{M}_\omega\right)^{p_{m+1}}\hspace{-.5mm}\matr{M}_{\omega\mu_m}\hspace{-1mm}\cdots \left(\matr{M}_{\mu_2}\right)^{p_{2}}\hspace{-.5mm}\matr{M}_{\mu_2\alpha}\hspace{-.5mm}\left(\matr{M}_\alpha\right)^{p_{1}}\hspace{-2mm},\hspace{-2mm}
\end{multline}
where the again loop numbers satisfy $\sum_{i=1}^{m+1}p_i=n-m$, with the restriction that the loop numbers on any loopless vertices are fixed to zero. Such a sequence is said to form a weak composition of $n-m$. 
We now note that for any weak composition of $n-m$, the following identity holds:
\begin{align}
\label{eq:fnExpDecomposed}
  \frac{1}{p_1!\,p_2!\cdots p_{m+1}!}\int_0^t \mathrm{d}t_m \cdots\int_0^{t_2}\mathrm{d}t_1\, (t-t_m)^{p_{m+1}}\cdots
  (t_2-t_1)^{p_{2}}\, t_1^{p_1}=\frac{t^n}{n!}.
\end{align}
This result -- which does not depend on the value of each individual $p_i$ -- is straightforwardly proved by induction on $m$. By substituting this identity into Eq.~\eqref{eqn:DecompOfMatrixExpI} and rearranging the order of summations we obtain
\begin{subequations}
\begin{align}
  \exp\big(t\,\matr{M}\big)_{\omega\alpha} &\hspace{-1mm}=\label{eqn:MatrixExpNestedConvolution1}\\
  &\hspace{-20mm}\sum_{m=0}^\infty\sum_{W_{\mathcal{G}_0;\alpha\omega;m}}\sum_{\{p_i\}=0}^\infty \int_0^t \hspace{-1.5mm}\mathrm{d}t_m\cdots\hspace{-1mm}\int_0^{t_2}\hspace{-1.5mm}\mathrm{d}t_1 \frac{\big[(t-t_m)\matr{M}_{\omega}\big]}{p_{m+1}!}^{\!p_{m+1}}\hspace{-5.5mm}\matr{M}_{\omega\mu_m}\hspace{-1mm}\cdots\frac{\big[(t_2-t_1)\matr{M}_{\mu_2}\big]}{p_{2}!}^{\!p_{2}}\hspace{-2mm}\matr{M}_{\mu_2\alpha} \frac{\big[t_1\matr{M}_{\alpha}\big]}{p_{1}!}^{\!p_{1}}\hspace{-2mm},\nonumber\\
&\hspace{-15mm}=\sum_{m=0}^{\infty}\sum_{W_{\mathcal{G}_0;\alpha\omega;m}}\hspace{-1mm}\int_0^t\hspace{-1.5mm} \mathrm{d}t_m\cdots\hspace{-1mm}\int_0^{t_2}\hspace{-1.5mm}\mathrm{d}t_1\,
 \mathrm{e}^{(t-t_m)\matr{M}_\omega}\matr{M}_{\omega\mu_m}\hspace{-.5mm}\cdots
  \mathrm{e}^{(t_2-t_1)\matr{M}_{\mu_2}}\,\matr{M}_{\mu_2\alpha}\,\mathrm{e}^{t_1\matr{M}_{\alpha}}.
\end{align}
\end{subequations}
This intermediate result proves Lemma \ref{thm3:ResummedMatrixExponential}. In order to continue, we note that this expression is an $m$-fold nested convolution.
To convert the convolution to a product we take the Laplace transform of both sides:
\begin{align}
 \hspace{-3mm} \mathcal{L}\big[\exp(t\,\matr{M})_{\omega\alpha}\big] &= \sum_{m=0}^\infty\sum_{W_{\mathcal{G}_0;\alpha\omega;m}}\mathcal{L}\big[f_\omega(t)\big]\,
  \matr{M}_{\omega\mu_m}\cdots \mathcal{L}\big[f_{\mu_2}(t)\big]\, \matr{M}_{\mu_2\alpha}\,\mathcal{L}\big[f_\alpha(t)\big], \\
 &\hspace{-10mm}= \sum_{m=0}^\infty\sum_{W_{\mathcal{G}_0;\alpha\omega;m}}\left[s\matr{I}-\matr{M}_\omega\right]^{-1}\,
  \matr{M}_{\omega\mu_m}\cdots\left[s\matr{I}-\matr{M}_{\mu_2}\right]^{-1}\, \matr{M}_{\mu_2\alpha}\,\left[s\matr{I}-\matr{M}_\alpha\right]^{-1},\label{eqn:MatrixExpSumOverLoopFreeWalks}
\end{align}
where the second line follows on noting the result $\mathcal{L}\big[\exp(t\matr{M}_\mu)\big] = \big[s\matr{I}-\matr{M}_\mu\big]^{-1}$.
As in the previous sections, we remark that Eq.~(\ref{eqn:MatrixExpSumOverLoopFreeWalks}) is a sum of walk contributions with an effective vertex weight of $c^{\textrm{eff}}_\mu=\left[s\matr{I}-\matr{M}_\mu\right]^{-1}$. On using Theorem \ref{PSresult} to turn Eq.~(\ref{eqn:MatrixExpSumOverLoopFreeWalks}) into a path-sum, we obtain
\begin{subequations}
\begin{align}
 &\hspace{-3mm} \mathcal{L}\big[\exp\big(t\,\matr{M}\big)_{\omega\alpha}\big] =\sum_{\ell=0}^{\ell_\mathrm{max}} \sum_{P_{\mathcal{G}_0;\alpha\omega;\ell}}\matr{F}_{\mathcal{G}\backslash\{\alpha,\ldots,\nu_\ell\}}[\omega]\,\matr{M}_{\omega\nu_\ell}\,\cdots\, \matr{F}_{\mathcal{G}\backslash\{\alpha\}}[\nu_2]\,\matr{M}_{\nu_2\alpha}\,\matr{F}_\mathcal{G}[\alpha]\,,\\
&\matr{F}_{\mathcal{G}}[\alpha] = \Big[s\matr{I}-\matr{M}_{\alpha} -
  \sum_{C_{\mathcal{G}_0;\alpha}}\matr{M}_{\alpha\mu_m} \matr{F}_{\mathcal{G}\backslash{\{\alpha,\ldots,\mu_{m-1}}\}}[\mu_m]\cdots \matr{F}_{\mathcal{G}\backslash\{\alpha\}}[\mu_2]\matr{M}_{\mu_2\alpha}\Big]^{-1},
\end{align}
\end{subequations}
where $\ell$ and $m$ are the length of the path and of the bare cycle, respectively. In this expression, the contribution from the loops is explicitly separated from that of the other bare cycles. 
\qquad\end{proof}

\subsection{\textbf{The matrix logarithm}}\label{subsec:MatrixLogProof}
\begin{proof}
In this section, we consider a matrix $\matr{M}\in\mathbb{C}^{D\times D}$ with no eigenvalue on the negative real axis. 
To prove Theorem \ref{thm:PathSumMatrixLog} we write the matrix logarithm as $
  \log \matr{M} = -\sum_{n=1}^\infty(\matr{I}-\matr{M})^n/n$. This series only converges if $\|\matr{I}-\matr{M}\| < 1$; nevertheless, the end result of this proof can be extended to matrices of arbitrary norm by analytic continuation. Theorem \ref{thm:PathSumMatrixLog} is therefore valid for all matrices regardless of norm. We introduce the auxiliary matrix $\bar{\matr{M}} \equiv \matr{I}-\matr{M}$ and rewrite the power series as
\begin{subequations}
\begin{align}
  \big(\log \matr{M}\big)_{\omega\alpha} &=-\sum_{n=1}^\infty \frac{1}{n}\sum_{W_{\mathcal{G}_0;\alpha\omega;n}} \bar{\matr{M}}_{\omega\nu_n}\cdots\bar{\matr{M}}_{\nu_3\nu_2}\bar{\matr{M}}_{\nu_2\alpha},
  \\ &\hspace{-16mm}=-\sum_{n=1}^\infty \frac{1}{n}\sum_{m=0}^n\sum_{W_{\mathcal{G}_0;\alpha\omega;m}}\sum_{\{p_i\}/ n-m} \hspace{-2mm}\big(\bar{\matr{M}}_\omega\big)^{p_{m+1}} \bar{\matr{M}}_{\omega\mu_m} \hspace{-1.5mm}\cdots \big(\bar{\matr{M}}_{\mu_2}\big)^{p_2} \bar{\matr{M}}_{\mu_2\alpha}\big(\bar{\matr{M}}_\alpha\big)^{p_1},\label{eqn:MatrixLogProof1}
\end{align}
\end{subequations}
with $\mathcal{G}$ the graph of the partition of $\bar{\matr{M}}$. The second equality is obtained by making the loops explicit through the same procedure as in \S\ref{subsec:MatrixPowerProof}.
Just as in the previous sections, the loop numbers form a weak composition of $n-m$. For any such composition, the following identity holds:
\begin{align}
\label{eq:fnLogDecomposed}
\frac{1}{n} = \int_0^1 x^{m-1}\, x^{p_1}x^{p_2} \,\cdots \, x^{p_{m+1}}\,\mathrm{d} x.
\end{align}
This identity allows the contributions from the infinite loop sums in Eq.~\eqref{eqn:MatrixLogProof1} to be evaluated in closed form. By restructuring the double summation we obtain
\begin{align}\label{eqn:MatrixLogProof2}
 &\big(\log \matr{M}\big)_{\omega\alpha} = -\delta_{\omega\alpha}\int_0^1 \bar{\matr{M}}_\alpha\left[\matr{I}-x\bar{\matr{M}}_\alpha\right]^{-1}\mathrm{d}x\\
 &-\sum_{m=1}^\infty\sum_{W_{\mathcal{G}_0;\alpha\omega;m}} \int_0^1 x^{m-1} \left[\matr{I}-x\bar{\matr{M}}_\omega\right]^{-1}\bar{\matr{M}}_{\omega\mu_m}\cdots \left[\matr{I}-x\bar{\matr{M}}_{\mu_2}\right]^{-1}\bar{\matr{M}}_{\mu_2\alpha}\left[\matr{I}-x\bar{\matr{M}}_{\alpha}\right]^{-1}\,\mathrm{d}x,\nonumber
\end{align}
where we have written $\left[\matr{I}-x\bar{\matr{M}}_\mu\right]^{-1}$ for $\sum_{p=0}^\infty (x \bar{\matr{M}}_\mu)^p$ and $\delta_{\omega\alpha}$ is a Kronecker delta. Note that the first part of this expression -- which represents contributions from walks consisting entirely of loops -- has a slightly different structure to the second part, owing to the absence of the term with zero loops when $\omega = \alpha$. Just as for the previously obtained matrix functions, the integrand of Eq.~(\ref{eqn:MatrixLogProof2}) is a sum of walk contributions with effective vertex weights $c^{\textrm{eff}}_\mu=\left[\matr{I}-x\bar{\matr{M}}_\mu\right]^{-1}$. On using Theorem \ref{PSresult}, Eq.~(\ref{eqn:MatrixLogProof2}) transforms to
\begin{subequations}
\begin{align}
  &\big(\log \matr{M}\big)_{\omega\alpha} =\\ &\hspace{2mm}\begin{dcases} \int_0^1\, \mathrm{d}x\,\, x^{-1}\, \left(\matr{I}-\matr{F}_{\mathcal{G}}[\alpha]\right), & \omega = \alpha, \\
  -\sum_{P_{\mathcal{G}_0;\alpha\omega}} \int_0^1\, \mathrm{d}x\,\, x^{\ell-1}\, \matr{F}_{\mathcal{G}\backslash\{\alpha,\ldots,\nu_\ell\}}[\omega] \bar{\matr{M}}_{\omega\nu_\ell} \cdots \matr{F}_{\mathcal{G}\backslash\{\alpha\}}[\nu_2]\, \bar{\matr{M}}_{\nu_2\alpha}\, \matr{F}_\mathcal{G}[\alpha],& \omega \neq \alpha,
\end{dcases}\nonumber
 \shortintertext{where $\ell$ the length of the path and}
&\matr{F}_{\mathcal{G}}[\alpha] \hspace{-.2mm}=\hspace{-1mm}\Big[\matr{I}- x\bar{\matr{M}}_\alpha-\hspace{-2mm}\sum_{C_{\mathcal{G}_0;\alpha}} x^m\,\bar{\matr{M}}_{\alpha\mu_m} \matr{F}_{\mathcal{G}\backslash\{\alpha,\ldots,\mu_{m-1}\}}[\mu_m]\,\cdots\, \bar{\matr{M}}_{\mu_3\mu_2}\matr{F}_{\mathcal{G}\backslash\{\alpha\}} [\mu_2]\bar{\matr{M}}_{\mu_2\alpha}\Big]^{-1}\hspace{-3mm},
\end{align}
\end{subequations}
with $m$ the length of the bare cycle. Note that we have again explicitly separated the contribution of the loops from that of the other bare cycles.
Theorem \ref{thm:PathSumMatrixLog} is now directly obtained upon replacing $\bar{\matr{M}}$ by $\matr{I}-\matr{M}$.
\qquad\end{proof}

\bibliographystyle{siam}

\end{document}